\documentclass[graybox]{svmult}
\usepackage{mathptmx}       
\usepackage{helvet}         
\usepackage{courier}        
\usepackage{type1cm}        
\usepackage{makeidx}         
\usepackage{graphicx}        
\usepackage{multicol}        
\usepackage[bottom]{footmisc}

\makeindex             

\usepackage{amsmath, amssymb, amscd, amsfonts} 

\def\qqqed{\hspace{\stretch{1}}  \qed}

\newcommand{\ep}[0]{\qqqed
\end{proof}
}

\newcommand{\NN}[0]{
{\mathbb{N}}
}

\newcommand{\RRd}[0]{
{\mathbb{R}^d}
}

\newcommand{\Hil}[0]{
\mathcal{H} }
\newcommand{\norm}[2]{
\left\| #2 \right\|_{#1} }

\newcommand{\BL}[1]{
{\mathfrak{B} \big( #1 \big)}
}

\newcommand{\AAa}[0]{
{\mathcal A} }

\newcommand{\HHmp}[0]{
{{\Hil}_w^p}
}
\newcommand{\HHmpdual}[0]{
{{\Hil}_{1/w}^q}
}
\newcommand{\HHmpdddd}[0]{
{\left( \HHmp \right)'}
}
\newcommand{\HHmi}[0]{
{{\Hil}_w^\infty}
}
\newcommand{\HHmidual}[0]{
{{\Hil}_{1/w}^1}
}

\newcommand{\llmp}[0]{
{\ell_w^p}
}

\newcommand{\Hn}[0]{
\Hil^{00}
}

\def\WWWfun{{\mathcal{W}}}
\newcommand{\WWW}[1]{
{\WWWfun \left( #1 \right)}
}

\def\VVVfun{{\mathcal{V} }}
\newcommand{\VVV}[1]{
{\VVVfun \left( #1 \right)}
}

\newcommand{\range}[1]{\mathsf{ran}\left( #1 \right)} 
\newcommand{\domain}[1]{\mathsf{dom}\left( #1 \right)} 
\newcommand{\identity}[1]{\mathsf{id}_{ #1 }}
\newcommand{\kernel}[1]{\mathsf{ker}\left( #1 \right)}

\def\Coosp{{\boldsymbol{\mathcal Co}}}
\newcommand{\Coo}[1]{{\Coosp \left( #1\right)}}

\def\Oosp{{\boldsymbol{\mathcal O}}}
\newcommand{\Oo}[1]{{\Oosp \left( #1\right)}}

\begin{document}

\title*{A Guide to Localized Frames and Applications to Galerkin-like Representations of Operators}
\titlerunning{Localized Frames and Galerkin} 
\author{Peter Balazs and Karlheinz Gr\"ochenig}
\institute{Peter Balazs \at Acoustics Research Institute, Austrian Academy of Sciences, Wohllebengasse 12-14, Vienna, 1040, Austria
\email{peter.balazs@oeaw.ac.at}
\and Karlheinz Gr\"ochenig \at Faculty of Mathematics, University of Vienna, Oskar-Morgenstern-Platz 1, Vienna, 1090, Austria \email{karlheinz.groechenig@univie.ac.at}
}
%
%
\maketitle

\abstract*{
This chapter offers a detailed survey on intrinsically localized frames and the corresponding matrix representation of operators. 
We re-investigate the properties of localized frames and the associated Banach spaces in full detail. 
We investigate the representation of operators using localized frames in a Galerkin-type scheme. 
We show how the boundedness and the invertibility of matrices and operators are linked and give some sufficient and necessary conditions for the boundedness of operators between the associated Banach spaces. 
}

\abstract{
This chapter offers a detailed survey on intrinsically localized frames and the corresponding matrix representation of operators. 
We re-investigate the properties of localized frames and the associated Banach spaces in full detail. 
We investigate the representation of operators using localized frames in a so-called Galerkin-type scheme. 
We show how the boundedness and the invertibility of matrices and operators are linked and give some sufficient and necessary conditions for the boundedness of operators between the associated Banach spaces. 
}
\section{Introduction}

Localized frames are ``good'' frames. More precisely, the concept of
localized frames was introduced in~\cite{gr04-1} in an attempt to understand
which properties render a frame useful. Whereas an abstract frame can
be viewed as a flexible coordinate system for a Hilbert space --- and only
for one Hilbert space! --- localized frames go beyond Hilbert spaces and yield
a description and characterization of a whole class of associated
Banach spaces. Indeed, the success of structured frames, such as Gabor
frames~\cite{feistro1}, wavelet frames~\cite{daubech1}, or frames of
translates~\cite{Cachka01}, is built on their capacity to describe modulation
spaces (Gabor), Besov-Triebel-Lizorkin spaces (wavelet), and their use
in sampling theory (translates). 

Gabor frames are used for the description and extraction of
time-frequency features of a signal. It would be a waste of possibilities to use them  to merely
determine  the $L^2$-membership of a
given function. Likewise wavelets are used to detect edges in an image
or compress an image, and not just for the expansion of a function in
$L^2$. In these applications one does not use arbitrary Gabor frames
or wavelet frames, but the window and the wavelet are usually
carefully designed so as to have some desired time-frequency
concentration or a small support and vanishing moments. Thus in such
applications the frames come with an additional property, namely some
form of \emph{localization}. 

The general theory of localized frames began as an attempt to
formulate an abstract concept of localization that can explain the
success of certain structured frames in applications~\cite{gr04-1}. Roughly
speaking, a frame $\Psi = \{ \psi _k : k\in K\}$ is called  localized, if its Gramian  matrix $G$ with
entries $G_{k,l}= \left< \psi_l, \psi_k \right>_{k,l \in K}$ possesses enough
off-diagonal decay. In the further developments of the concept
powerful methods of Banach algebra theory were used, and nowadays, and
also in this survey, we call a frame localized, if its Gramian belongs
to a solid, inverse-closed Banach algebra of matrices~\cite{forngroech1}. \\

Localized frames possess many properties   that are not available for
general frames. 

(i) To every localized frame can be associated a class of Banach
spaces, the so-called coorbit spaces \cite{fegr89,fegr89-1}. Roughly speaking, the coorbit
space $\mathcal{H}^p_w$ contains all elements $f$ such that the
sequence $\langle f , \psi _k\rangle _k$ belongs to the weighted $\ell
^p$-space. For the standard structured frames  one obtains well-known
and classical families of functions spaces, namely, the 
modulation spaces are associated to Gabor frames~\cite{gr01}, and the Besov spaces
are associated to wavelet frames~\cite{meyer1,raul11}. In this chapter we will
explain the construction of the coorbit spaces and derive many of
their properties.  

(ii) Localized frames possess nice dual frames. Technically speaking,
the canonical dual frame possesses the same localization. In fact,
this is a fundamental statement about localized frames,  and the
method of proof (matrix algebras, spectral invariance) 
has motivated intensive research in  Banach algebra theory~\cite{grle06,grrz08,grkl10,gr10-2}.

(iii) Elements in coorbit spaces associated to a localized frame
possess good approximation properties~\cite{corgroech04}.  In fact,
the results on  nonlinear $N$-term approximations and on the fast
convergence of the iterative frame algorithms are based fundamentally
on the localization property and do not hold for arbitrary frames. 

(iv) Localized frames often possess a characterization ``without
inequalities'' ~\cite{gro07ineq,grorro15,degrro16}. These results have lead to
strong results about Gabor frames that have been out of reach with
conventional methods. 

(v) Every localized frame is a finite union of Riesz
sequences~\cite{gr03-4}. This is a 
special case of the Feichtinger conjecture and  was proved early on
with an easy proof,  whereas the recent  full proof of the Feichtinger conjecture is a
monumental result in functional analysis and occupies a completely
different mathematical universe~\cite{maspsr15}. 

(vi) General frames may be used to describe and discretize operators
and operator equations, 
and thus have led to an important line of frame research in numerical
analysis. In the so-called Galerkin approach an operator $O$
is discretized with respect to a frame by associating a matrix $M$,
with  
$M_{k,l} = \left< O \psi_k, \phi_l \right>$  with respect to given
frames $\Psi , \Phi $. The standard  discretization  uses  bases, but recently
also frames have been used \cite{xxlframoper1}. 
The Galerkin approach works particularly well 
when the corresponding matrix is sparse. The most famous example is
the sparsification of singular integral operators by means of wavelet
frames (or bases)~\cite{beycoif91,me90-8}. This work  has led to many adaptive methods
\cite{harbr08,Stevenson03,Dahlkeetal07c}. 
In this regard the time-frequency analysis of pseudodifferential
operators by means of Gabor frames is a particularly
successful example  of the application of  localized frames: certain symbol classes containing the
H\"ormander class $S^0_{0,0}$ can be completely characterized by the
off-diagonal decay of the associated matrix~\cite{grrz08}.
Subsequently Gabor frames were pivotal for the sparsification of
Fourier integral operators and certain Schr\"odinger propagators in
the work of the Torino group~\cite{cogrniro13-1,coniro09,coniro15,coniro15-2}
On a more abstract level, localized frames have been heavily used in
the adaptive frame methods for the solution of operator
equations in~\cite{dafora05,Dahlkeetal10}. \\

This chapter will offer a survey of localized frames. Of the many
possible topics we will focus on the theory of the associated coorbit
spaces and on the Galerkin discretization of operators with respect to
localized frames. 
 We will mainly
explain the abstract theory and focus on the formalism of localized
frames. These  aspects have not received as much attention as  other
topics and justify a more thorough treatment. Several results may even
claim some novelty, for instance, the inclusion relations
of coorbit spaces and the explicit relation between the mapping
properties of an operator and of its corresponding matrix seem to be
new. 

Although the topic of localized frames is eminently applied, we will
choose a formalistic point of view and develop and explain the
\emph{formalism} of localized frames, their coorbit spaces, and the Galerkin
discretization.

The motivation for this formal approach, and for this entire chapter,
comes from different readings of the original sources~\cite{gr04-1,forngroech1} and many
ensuing discussions between the authors. One of us (K.\ G.) claimed
that ``this is folklore and known'', while the other (P.\ B.) would
point  out --- and rightly so ---  that the results he needed and wanted to understand
in detail were not formulated in the publications. P.\ B.\ strongly
argued that he needed a general abstract formalism in order move on
to the real applications in acoustic applications as e.g. in~\cite{riecketal11}. 
The result of our discussions is this summary of localized frames
with its emphasis on the formalism. We hope that this point of view
will also benefit other readers and users of localized frames. \\

This chapter  is organized as follows. In Section \ref{sec:prelnot0}
we collect some preliminary definitions and notation and then
introduce the concept of localization frames. %
Section  \ref{sec:localfr0} 
is devoted to the study of the associated
coorbit spaces and the canonical operators associated to a frame. 
In Section \ref{sec:matreplocalfr0} we describe the Galerkin
approach and discuss the formalism 
of  matrix representations of operators with respect to localized
frames. We end with a short outlook in Section \ref{sec:soopeq0}.

\section{Preliminaries and Notation} \label{sec:prelnot0}

For a standard reference to functional analysis and operator theory refer e.g. to \cite{conw1}. 
 We denote by $\BL{X}$ the Banach algebra of  bounded operators on the normed space $X$.
 We will write $\norm{X \rightarrow Y}{T}$ for the operator norm of a bounded operator $T: X \rightarrow Y$, or just $\norm{}{T}$, if the spaces are clear. 
 We will use the same notation for the inner product of a Hilbert space $\left< . , . \right>_\Hil$ and for the duality relation of two dual sets $B,B'$, $\left< . , . \right>_{B,B'}$. 
 If there is no chance of confusion, we will just use the notation $\left< . , . \right>$ for that. 

Let $A \in \BL{\Hil_1,\Hil_2}$ with closed range. Then there exists a unique bounded operator $A^\dagger : \Hil_2 \rightarrow \Hil_1$ for which
$ A A^\dagger f = f, \forall f \in ran(A) $ and $\kernel{A^\dagger} = \left( \range{A} \right)^\bot$. This $A^\dagger$ is called the (Moore-Penrose) \emph{pseudoinverse} of $A$. See e.g. \cite{olepinv}.\\

\subsection{Sequence spaces}

We use the canonical notation of sequence spaces $\ell^p = \ell^p (K)$
consisting of  sequences on a  countable and separable  index set $K \subseteq \RRd$.
By an abuse of notation, but for greater consistency, we define $\ell^0$ as those sequences for which $\lim \limits_{k\rightarrow \infty} c_k = 0$. (Usually this space is denoted by $c_0$.)
We denote the set of sequence with only finitely many non-zero entries by $\ell^{00}$ (usually denoted by $c_{00}$).

A weight is a (strictly) positive sequence $w = (w_k)_{k \in K}$, $w_k > 0$. Then we define the weighted space $\ell_w^p$ by
$c \in \ell_{w}^{p} \Longleftrightarrow w \cdot c \in \ell^p$ with norm \mbox{$\norm{\ell_{w}^{p}}{c} = \norm{\ell^p}{c \cdot w}$.}
So for the 
weighted sequence spaces and $1 < p \le 2$ and $2 \le q < \infty$ we get
\begin{equation} \label{sec:scalseq1}   \ell^{00} \subseteq \underbrace{ \ell_w^1 \subseteq \ell_w^p \subseteq \ell_w^2 \subseteq \ell_w^q \subseteq \ell_w^0}_{(*)} \subseteq \ell_w^{\infty}
\end{equation} 
where the middle $(*)$ is a chain of dense Banach spaces. 
$ \ell^{00}$ is dense in all $\ell_w^p$ for $1 \le p < \infty$ and $p = 0$, and weak-* dense in $\ell_w^\infty$. 
Clearly $\ell_w^0 = \overline{ \ell^{00}}^{\norm{\ell_w^\infty}{.}}$. 

For $1 \le p < \infty$ and $1 = \frac{1}{p} +\frac{1}{q}$ we know that $\left( \ell_w^p \right)' \cong \ell_{1/w}^q$  with the duality relation 
\begin{equation}\label{sec:dualqeightseq1}\left< c_k , d_k
  \right>_{\ell_w^p,\ell _{1/w}^q} = \left< w_k c_k , \frac{1}{w_k} d_k \right>_{l^p,l^q} = \sum \limits_k c_k \overline{d_k}.
\end{equation}   
For $p = \infty$ this is only true in the K{\"o}the dual sense \cite{Koethe1}. We also have $\left(l_w^0\right)' \cong \ell_{1/w}^1 $.

\subsection{Frames}

A sequence $\Psi = \left( \psi_k\right)_{k\in K}$ in a separable Hilbert space $\mathcal{H}$ is a {\sl frame} for $\mathcal{H}$, if there exist positive
constants $A$ and $B$ (called lower and upper frame bound, respectively) that satisfy
\vspace{-2mm}
\begin{equation} \label{sec:frambasdef1}
A \|f\|^2\leq\sum_{k\in K}|\langle f,\psi_k\rangle|^2 \leq
B \|f\|^2\,\,\,\,\,\forall f \in \mathcal{H} .
\vspace{-2mm}
\end{equation}
A frame where the two bounds can be chosen to be equal, i.e.,  $A = B$, is called {\em tight}.
 In the following we will denote by $\Psi =(\psi_k)$ and $\Phi = (\phi_k)$ the corresponding sequences in $\Hil$. 

By $C_\Psi : \Hil \rightarrow \ell^2$ 
 we denote the {\em analysis operator} defined by $\left(C_\Psi f\right)_k = \left< f , \psi_k\right>$.
The adjoint of $C_\Psi$ is the {\em synthesis operator} $D_\Psi (c_k) = \sum_k c_k \psi_k$.
The {\em frame operator} $S_\Psi = D_\Psi C_\Psi$ can be written as $S_\Psi f = \sum_k \left< f , \psi_k\right> \psi_k$, it is positive and invertible\footnote{Note that those 'frame-related' operators can be defined as possibly unbounded operators for any sequence \cite{xxlstoeant11}.}.
By using the {\em canonical dual frame} $(\tilde \psi_k)$, $\tilde
\psi_k = S^{-1} \psi_k$ for all $k$, we obtain the  reconstruction
formula
 $$f = \sum \limits_k \left< f , \psi_k\right> \tilde \psi_k = \sum
 \limits_k \left< f , \tilde \psi_k\right> \psi_k \qquad  \text{ for
   all }  f \in
 \Hil \, .$$
 Any sequence $\Psi _d =  (\psi_k^d)$ for which $C_{\Psi ^d}$ is
 bounded on $\mathcal{H}$ and where such a reconstruction holds is called a
 dual frame. 

The {\em Gram matrix} $G_{\Psi}$ is defined by $\left(G_{\Psi}\right)_{k,l} = \left< \psi_l, \psi_k \right>$. 
This matrix defines an operator on $\ell^2$ by matrix multiplication, corresponding to $G_{\Psi} = C_\Psi D_\Psi$. Similarily we can define the {\em cross-Gram matrix} ${\left(G_{\Psi, \Phi}\right)}_{k,l} = \left< \phi_l ,\psi_k \right>$ for two frames $\Phi$ and $\Psi$.
Clearly 
$$G_{\Psi, \Phi} c = \sum \limits_l \left(G_{\Psi,\Phi}\right)_{k,l}
c_l = \left< \sum \limits_l c_l \phi_l , \psi_k \right> = C_\Psi
D_\Phi c \, .$$

If, for the sequence $\Psi$,  there exist constants $A$, $B >0$ such that the inequalities
$$ A \norm{2}{c}^2 \le \norm{\Hil}{\sum \limits_{k \in K} c_k \psi_k}^2 \le B \norm{2}{c}^2 $$
are  fulfilled, $\Psi$ is called a {\em Riesz sequence}. If $\Psi$ is complete, it is called a {\em Riesz basis}.

\subsubsection{Banach frames} \label{sec:banchfram0}

The concept of frames can be extended to Banach spaces \cite{gr91,chst03,Casazza2005710}: 

Let $X$ be a Banach space and $X_d$ be a Banach space of scalar sequences.
A sequence
 $(\psi_k)$ in the dual $X^\prime $ is called a $X_d$-frame
for the Banach space $X$  if there exist constants $A, B>0$
such that
\begin{equation} \label{sec:banframbasdef1}
A\|f\|_{X}\leq \norm{X_d}{ \langle f, \psi _k\rangle _{k\in K}} \leq
B\|f\|_{X}\quad\textrm{for all}\quad
 f \in X.
\end{equation}

An $X_d$-frame is called a Banach frame with respect to a sequence space $X_d$, if there exists a bounded reconstruction operator $R : X_d \rightarrow X$, such that 
$ R \left( \psi_k(f) \right) = f$ for all $f \in X$.
If  $X_d = \ell^p$ for $1\le p \le\infty$, we speak of  $p$-frames, respectively $p$-Banach frames.
 The distinction between $X_d$-frames and Banach frames will disappear for localized frames. The norm equivalence \eqref{sec:banframbasdef1} already implies the existence of a reconstruction operator for $X$, in this setting.

\subsubsection{Gelfand triples}

Let $X$ be a Banach space and $\Hil$ a  Hilbert space. Then the triple $(X, \Hil, X')$ is called a Banach Gelfand triple, if
$X \subseteq \Hil \subseteq X'$, where  $X$ is dense in $\Hil$,  and $\Hil$ is $w^*$-dense in $X'$.  
The prototype of such a triple is $(\ell^1, \ell^2, \ell^\infty)$.

A frame for $\Hil$ is called a {\em Gelfand frame} \cite{dafora05} for
this triple if there exists a Gelfand triple of sequence spaces $(X_d,
\ell^2, X_d')$, such that the synthesis operator $D_{\Psi} : X_d
\rightarrow X$ and the analysis operator $C_{\tilde \Psi} : X
\rightarrow X_d$ are bounded.

Now for a Gelfand frame $\Psi$ for the Gelfand triple $(X, \Hil, X')$ with the sequence spaces $(\ell^1, \ell^2, \ell^\infty)$, we define the coorbit space $\Coo{\ell^p, \Psi} = \left\{ f \in X' :  C_\Psi f \in \ell^p  \right\}$. Similarly, one could define the orbit spaces $\Oo{\ell^p, \Psi} = \left\{ D_\Psi c \text{ for } c \in \ell^p \right\}$. We refer to \cite{fegr89} for an early example and the terminology of coorbit spaces.

\section{Localization of frames} \label{sec:localfr0}
In this section we introduce the concept of localized frames and define the corresponding family of coorbit spaces.
In Subsection \ref{sec:hhwwinf0} we treat the maximal space $\HHmi$ in detail. 
In Subsection \ref{sec:dual0} we show the duality relations of these spaces . 
In Subsection \ref{sec:framrelop0} we study the frame-related operators.
\\

We call a Banach *-algebra $\AAa$ of infinite matrices (over the index set $K$) a {\em solid spectral matrix algebra}, if
\begin{enumerate}
\item[(i)] $\AAa \subseteq \BL{\ell^2}$.  
\item[(ii)] $\AAa$ is inverse-closed in $\BL{\ell^2}$, i.e.,  $A \in \AAa$
  and $A$ is invertible on $\ell^2$, then $A^{-1} \in \AAa$.  
\item[(iii)] $\AAa$ is {\em solid}, i.e.,  $A \in \AAa$ and $| b_{k,l} | \le
  |a_{k,l}|$, then $B = (b_{k,l}) \in \AAa$ and $\norm{\AAa}{B} \le
  \norm{\AAa}{A}$. 
\end{enumerate}
Several examples, e.g., the Jaffard class or a Schur-type class, can be found in \cite{forngroech1}.  
In these examples localization is defined by some off-diagonal decay
of the Gram matrix. 
For the systematic construction of spectral matrix algebras we refer
to ~\cite{grle06,grrz08,grkl10,Sun07a}, a survey on spectral invariance is
contained in \cite{gr10-2}.   

\begin{definition} \label{defloc} 
Let $\AAa$ be solid spectral matrix algebra. 
Two frames $\Psi$ and $\Phi$ are called $\AAa$-localized with respect to each other, if their cross-Gram matrix $G_{\Psi,\Phi}$ belongs to $\AAa$. 
If $G_{\Psi, \Phi} \in \AAa$, we write $\Psi \sim_{\AAa} \Phi$.

A single frame 
$\Psi = (\psi_k)$ is called (intrinsically) $\AAa$-localized, if $\Psi \sim_{\AAa} \Psi$.  
\end{definition}

Alternative definitions of localized frames can be found in
\cite{fora05} (continuous frames),   \cite{gr03-3,gr04-1}
(localization with respect to a Riesz basis), \cite{bacahela06} ($\ell^p$-self-localization)  or
\cite{gepe11,pe15-2,fefupe15} (localization in terms of the intrinsic
metric on a manifold). Although all these concepts have their merits,
we will focus on the intrinsic localization of Definition~\ref{defloc}. 

The following connection holds for any chosen dual frame $\Phi^{d}$ \cite{forngroech1}:
\begin{equation} \label{sec:nearequiv0}
\Psi \sim_{\AAa} \Phi, \Phi^{d} \sim_{\AAa} \Xi \Longrightarrow \Psi \sim_{\AAa} \Xi.
\end{equation}

A weight $w$ is called $\AAa$-admissible, if every $A \in \AAa$ can be
extended to a bounded operator on $\ell_w^p$ for all 
$1 \le p \le \infty$ i.e,
$\AAa \subseteq \bigcap \limits_{1 \le p \le \infty} \BL{\ell_w^p}$. 

In the following, $\AAa$ is always a solid spectral Banach algebra of matrices on $K$.
Since $\AAa$ is a Banach *-algebra, if $w$ is $\AAa$-admissible, then $1/w$ is admissible, too. This is because for $A : \ell_w^p \rightarrow \ell_w^p$, we have $A^*: \ell_{1/w}^q \rightarrow \ell_{1/w}^q$ for $q > 1$. For $q = 1$, this argument is valid using the pre-dual.\\[-3mm]

\begin{definition}
Let $\Hn = \left\{ f = \sum \limits_{k} c_k \psi_k : c \in  \ell^{00} \right\}$ be the subspace of all finite linear combinations over $\Psi$. 

For $1 \le p < \infty$ define $\HHmp(\Psi,\tilde \Psi)$ as the completion of $\Hn$  with respect to the norm 
$$\norm{\Hil_w^{p}}{f} = \norm{\ell ^p_w}{C_{\widetilde \psi}(f)} .$$ 

Let $\Hil_w^{0}$ be the completion of $\Hn$ with respect to the norm
$$\norm{\Hil_w^{0}}{f} = \norm{\Hil_w^{\infty}}{f} = \norm{\ell^\infty_w}{C_{\tilde \psi}(f)}.$$

In Section~\ref{sec:hhwwinf0} we will define  the space $\Hil_w^{\infty}$  as a weak$^*$ completion
with respect to the metric  $\norm{\ell^\infty_w}{C_{\tilde
    \psi}(f)}$. Alternatively, we may define it as the bidual $\Hil
_w^\infty = (\Hil_w^{0})^{**}$. 
\end{definition}
We note right away that $\Hil ^p _w \subseteq \Hil ^q_w \subseteq \Hil
^0_w$ for $1\leq p \leq q$. 

As a consequence of this definition the analysis operator can be
extended to a bounded operator from $\Hil_w^p$ into $\ell_w^p$. 
\\

The main results in \cite{forngroech1} are summarized below. The first one describes the independence of $\Hil_w^p ( \Psi, \tilde \Psi)$ of the defining frame $\Psi$. 
\begin{proposition}[\cite{forngroech1}] \label{sec:equivspac0}  Let $\Phi$ and $\Psi$ be frames for $\Hil$ and $\Phi^d$ and $\Psi^d$ dual frames. 
If $\Psi^d \sim_{\AAa} \Psi$, $\Phi^d \sim_{\AAa} \Psi$ and $\Psi^d \sim_{\AAa} \Phi$, then $\Hil_w^p ( \Psi, \Psi^d) = \Hil_w^p ( \Phi, \Phi^d)$ with equivalent norms for all $1 \le p \le \infty$. 
\end{proposition} 
The proof of this result relies on the algebra properties of $\AAa$ and identities for Gram matrices. By this result 
we may therefore write unambiguously $\HHmp :=  \Hil_w^p ( \Psi, \tilde \Psi) = \Hil_w^p ( \tilde \Psi, \Psi)$. \\

In particular, let $\Phi = \Psi$. For a frame $\Psi$ and its dual $\Psi^d$, which are $\AAa$-localized with respect to each other, it can be shown that they are automatically Banach frames for all involved associated Banach spaces:
\begin{theorem}[\cite{forngroech1}] \label{sec:reconst0} 
Assume that $\Psi \sim_{\AAa} \Psi^d$. Then $\Psi$ is a Banach frame for $\Hil_w^p ( \Psi, \Psi^d)$ for $1 \le p < \infty$ or $p = 0$. The reconstructions
 $f = \sum \limits_{n \in N} \left< f , \psi_n \right> \psi_n^d$ and  $f = \sum \limits_{n \in N} \left< f , \psi_n^d \right> \psi_n$ converge unconditionally for $1 \le p < \infty$.
\end{theorem}

The assumptions of Proposition \ref{sec:equivspac0} can be weakened for the canonical dual frame, because it can be shown that an intrinsically localized frame is automatically localized with respect to its canonical dual. 
As a consequence an intrinsically localized frame is automatically a Banach frame for all associated Banach spaces. 
This is the main theorem about localized frames:
\begin{theorem}[\cite{forngroech1}] \label{sec:origrecon1} Let $\Psi$ be an intrinsically $\AAa$-localized frame, $\Psi \sim_{\AAa} \Psi$. Then 
$$\tilde \Psi \sim_{\AAa} \tilde \Psi \text{ and } \Psi \sim_{\AAa} \tilde \Psi .$$
As a consequence, $\Hil_w^p ( \Psi, \tilde \Psi) = \Hil_w^p ( \tilde \Psi, \Psi)$ and $\Psi$ is a $p$-Banach frame for $\Hil_w^p ( \Psi, \tilde \Psi)$ for $1 \le p < \infty$ or $p = 0$. The reconstructions
\begin{equation} \label{sec:reconform1}
f = \sum \limits_{n \in N} \left< f , \psi_n \right> \tilde \psi_n \text{ and } f = \sum \limits_{n \in N} \left< f , \tilde  \psi_n \right> \psi_n
\end{equation} converge unconditionally in $\HHmp$ for $1 \le p < \infty$.
\end{theorem}

Therefore the norms $\norm{\llmp}{C_\Psi f}$ and $\norm{\llmp}{C_{\widetilde \Psi} f}$ are equivalent, and the inequalities 
\begin{equation} \label{sec:syngrammop2}
\frac{1}{\norm{\llmp \rightarrow \llmp}{G_{\widetilde \Psi}}} \norm{\HHmp}{f} \le \norm{\llmp}{C_\Psi f} \le  \norm{\llmp \rightarrow \llmp}{G_{\Psi}} \norm{\HHmp}{f}.
\end{equation}
are valid for $1 \le p < \infty$ and $p = 0$

The unconditional convergence of the reconstruction formula \eqref{sec:reconform1} implies that both synthesis operators $D_\Psi$ and $D_{\widetilde \Psi}$ map $\ell^p_w$ onto $\HHmp$ for $1 \le p < \infty$ and $p = 0$. 
\sloppy Consequently, an equivalent norm on $\HHmp$ is given  by 
$$\inf \left\{ \norm{\ell_w^p}{c} : f = D_{\psi} c \right\} \qquad
\text{ for } f \in \HHmp \, .$$ 

In particular this means that the orbit 
 and co-orbit 
 definitions of $\HHmp$ coincide. 
\\

The best studied  examples of intrinsically localized frames are the
following.
\begin{enumerate}
\item[(i)] Frames of translates \cite{Cachka01,corgroech04},
\item[(ii)] Frames of reproducing kernels in a shift-invariant space~\cite{gr04-1,sun10a}
\item[(iii)] Gabor frames~\cite{feistro1,gr04-1,forngroech1},
\item[(iv)] Frames of reproducing kernels in (general) Bargmann-Fock
  spaces~\cite{Lin01},
\item[(v)] Wavelet frames that are orthogonal across different
  scales~\cite{corgroech04,fland99}. 
\end{enumerate}

However, not all useful frames are localized in the
sense of Definition~\ref{defloc}, among them are general wavelet frames,
frames of curvelets, shearlets, frames on manifolds, etc. Although
these frames do possess some form of localization, they are not part
of our theory of localized frames. 
While many of the constructions
discussed in this chapter, such as the definition and characterization
of coorbit spaces, can be carried out by hand or with different
techniques,   the main results for localized frames are not available for
wavelets or  curvelets. For  instance,   the decisive Theorem~\ref{sec:origrecon1}
and most of its consequences are  false for wavelet frames and their
many generalizations. 

\subsection{$\Hil_w^\infty$ as a normed space} \label{sec:hhwwinf0}

In the following we will focus on the theory of the coorbit spaces $\HHmp \big( \Psi, \widetilde \Psi \big)$ that are associated to a localized frame. We start with the ``distribution space'' $\HHmi$ and offer a thorough treatment. 
In \cite{forngroech1} ``the rigorous discussion was omitted to avoid tedious technicalities.''

Let $w$ be an $\AAa$-admissible weight. We define $\Hil_w^\infty$  as
a certain weak$^*$ completion of $\Hil$. 
We say that  two sequences $(f_n)$ and $(g_n)$ in $\Hil $  are equivalent, denoted by $f_n \sim g_n$, 
if $\langle f_n - g_n, \tilde \psi _k \rangle \to 0$ as $n\to \infty
$. Alternatively, $f_n - g_n \to 0$ in the $\sigma (\Hil, \Hil
^{00})$-topology.
 
\begin{definition} \label{sec:defHHmi1}
 We define    $\HHmi$ as the set of  equivalence classes of
sequences  $f = \left[ f_n \right]$, such that 
\begin{enumerate}
\item[(i)] $f_n \in \Hil$ for all $n \in \NN$,  
\item[(ii)] $\lim _{n\to \infty }   \left< f_n, \tilde \psi _k
  \right> = \alpha _k $ exists for all $k\in K$, 
\item[(iii)] $\sup \limits_n \norm{\ell_w^\infty}{C_{\widetilde \Psi} f_n} < \infty$.
\end{enumerate}
\end{definition}

 In this way $\Hil ^\infty _w$ is well-defined. The definition of $f$
 is independent of its representative. Indeed, if $f= [f_n]$ and
 $f_n\sim g_n$, then $\alpha _k = \lim _{n\to \infty } \langle f_n , \tilde \psi
 _k \rangle = \lim _{n\to \infty } \langle g_n , \tilde \psi
 _k \rangle$. 

Furthermore,  condition (iii)  implies that 
$| \langle f_n, \tilde \psi _k \rangle | w_k \leq C$ for all $n\in
\mathbf{N}$ and $k\in K$, consequently, $ |\alpha _k| w_k  = \lim _{n\to \infty }   |\left< f_n, \tilde \psi _k
  \right>| w_k   \leq C$ and thus $\alpha \in \ell
  ^\infty _w$. 
Now, write $\langle f , \tilde \psi _k \rangle = \alpha _k$.  and set
 \begin{equation}
   \label{eq:c35}
   \|f \|_{\Hil ^\infty _w} = \|\alpha \|_{\ell _w^\infty} \, . 
 \end{equation}
Therefore $C_{\widetilde \Psi} : \Hil_w^\infty \rightarrow \ell_w^\infty$ is a bounded operator. 

Clearly, \eqref{eq:c35} defines a  seminorm, because limits are linear and $\|\cdot
\|_{\ell ^\infty _w}$ is a norm. Now assume that 
$\norm{\Hil_w^\infty}{f}= 0$. This means that for every representative
$[f_n]$ of $f$ we have $\lim _n \langle f_n, \tilde \psi _k\rangle = 0$, or
equivalently  $f_n \sim 0$. 
Thus $f = 0 $  in $\Hil$ and $\norm{\Hil_w^\infty}{ \cdot}$
is indeed a norm. 

\begin{lemma}\label{completec} \begin{itemize}
\item[(i)] The map $f= [f_n] \mapsto      \alpha  = (\lim _n \langle f_n, \tilde \psi
_k \rangle)_{k\in K}$ is an isometric isomorphism from $\Hil ^\infty _w$ onto the subspace 
$V_\Psi = \{\alpha \in \ell ^\infty _w : \alpha = G_{\tilde \Psi,
  \Psi )} \alpha \} $.
\item[(ii)] The subspace $V_\Psi$ is closed in $\ell ^\infty _w$ and thus
$\Hil ^\infty _w$ is complete.  
\end{itemize}
 \end{lemma}
 \begin{proof}
(i) For $f\in \Hil $ we interpret the reconstruction formula $f = \sum
\limits_l \langle f, \tilde \psi _l \rangle \psi _l$ weakly as $\sum \limits_l \langle f,
\tilde \psi _l \rangle \, \langle \psi _l , \tilde \psi _k\rangle $,
or in operator notation as 
\begin{equation}
  \label{eq:c37}
 C_{\tilde \Psi } f = G_{\tilde \Psi , \Psi } C_{\tilde \Psi } f   \, .
\end{equation}
Now let $f = [ f_n] \in  \Hil ^\infty _w $ as in Definition \ref{sec:defHHmi1}. This is a sequence of vectors
$f_n \in \Hil$ such
that 
 $\lim _{n\to \infty }   \left< f_n, \tilde \psi _k
  \right> = \alpha _k$  and $\|C_{\tilde \Psi }
f_n \|_{\ell ^\infty _w} \leq C$ for all $n$.  This means that the
  sequence $C_{\tilde \Psi } f_n$ converges pointwise to $\alpha $ and
  is dominated by the sequence  $ C
  (w_l^{-1})_{l}\in \ell ^\infty _w$.    By
  dominated convergence it now follows that (again with pointwise
  convergence) 
$$
\alpha = \lim _{n\to \infty } C_{\widetilde \Psi } f_n = \lim _{n\to
  \infty } G_{\widetilde \Psi , \Psi } C_{\tilde \Psi } f_n =
	G_{\widetilde
  \Psi , \Psi } \lim_{n\to \infty } C_{\widetilde \Psi } f_n = 
	G_{\widetilde
 \Psi , \Psi } \alpha \, .
$$
Consequently, the limiting sequence $\alpha \in \ell ^\infty _w$
satisfies $\alpha = G_{\tilde \Psi , \Psi }\alpha $ and $\alpha \in
V_\Psi $.

 Conversely, let $\alpha \in V_\Psi $. Choose a sequence $F_n$ of
 finite subsets of $K$, such that $F_n \subset F_{n+1}$ and $\bigcup
 _{n=1}^\infty F_n = K$ and define 
$$
f_n = \sum _{l \in F_n } \alpha _l \psi _l  \in \Hil \, .
$$
Then clearly 
$$
\lim _{n\to \infty } \langle f_n, \tilde \psi _k\rangle = \lim _{n\to
  \infty } \sum _{l\in F_n} \alpha _l \psi _l  = G_{\tilde \Psi , \Psi
} \alpha = \alpha \, , 
$$
and $\sup _{k} | \langle f_n , \tilde \psi _k \rangle | w_k \leq C
\|\alpha \|_{\ell ^\infty _w}$. This means that $f= [f_n ]\in \Hil
  ^\infty _w$, and as a consequence the map $[f_n] \in \Hil ^\infty _w
  \mapsto \alpha \in V_\Psi $ is an isometric isomorphism. 

(ii) Assume that $\alpha _n  \in V_\Psi $ and $\alpha \in \ell ^\infty
_w$ such that $\| \alpha _n - \alpha \|_{\ell ^\infty _w} \to
0$. Since $ G_{\tilde \Psi , \Psi
} $ is bounded on $\ell ^\infty _w$,  we obtain that 
$ \alpha = \lim _{n\to \infty } \alpha _n = \lim _{n}  G_{\tilde \Psi , \Psi
} \alpha _n =  G_{\tilde \Psi , \Psi
} \alpha $, whence $\alpha \in V_\Psi $ and $V_\Psi $ is a
(norm)-closed subspace of $\ell ^\infty _w$. By the identification
proved in (i), $\Hil ^\infty _w$ is therefore complete.  
\ep 

Switching the roles of $\Psi$ and  $\widetilde \Psi$ in Definition
\ref{sec:defHHmi1}, we see that  $C_\Psi$ is an isometry between
$\HHmp(\widetilde \Psi, \Psi)$ and a closed subspace of $\ell_w^\infty$. By Equation \eqref{sec:syngrammop2} the corresponding norms are equivalent and so $C_\Psi$ is bounded from $\HHmp$ into $\ell_w^\infty$.

We next  verify the unconditional weak $^*$-convergence of the
sum  $D_\Psi c = \sum \limits_{k \in K} c_k \psi_k$ for $(c_k) \in \ell_w^\infty$. Let $\epsilon >0$
and $g\in \Hil ^{00}$. Choose a finite set $H_0$, such that 
$\sum _{k \not \in H_0} | \langle g, \psi _k\rangle | w_k ^{-1} <
\epsilon / \|c\|_{\ell ^\infty _w}$. Now let $H, J\subseteq K$ be two
finite sets such that $H\supseteq H_0$ and $J\supseteq H_0$. Then
$H\setminus J \cup J \setminus H \subseteq K \setminus H_0$, and
therefore 
\begin{align*}
 | \langle \sum \limits_{k \in J} c_k \psi_k &- \sum \limits_{k \in H}
 c_k \psi_k, g \rangle | = | \langle \sum \limits_{k \in 
  H\setminus J \cup J \setminus H} c_k \psi_k, g \rangle | \\
& \leq \|c\|_{\ell ^\infty _w} \,  \sum \limits_{k \not \in H_0} |
\langle \psi_k, g \rangle | w_k^{-1}  < \epsilon \, .  
\end{align*}
Thus  the series for $D_\Psi c$ converges weak-$^*$
unconditionally.  
Furthermore, since $G_{\tilde \Psi , \Psi }$ is bounded on $\ell
^\infty _w$ by the basic assumption on $\Psi $ and Theorem~\ref{sec:origrecon1}, we
also deduce the boundedness of $D_\Psi $ as follows: 
$$ \norm{\Hil_w^\infty}{D_\Psi c} = \norm{\ell_w^\infty}{C_{\tilde
    \Psi} D_\Psi c} = \norm{\ell_w^\infty}{G_{{\tilde \Psi},\Psi }c}
\leq \|G_{{\tilde \Psi},\Psi } \|_{\ell ^\infty _w \to \ell ^\infty _w
} \|c\|_{\ell ^\infty _w} < \infty. $$ 

The following lemma summarizes the properties of $\HHmi$. 
\begin{lemma} \label{sec:propinft1} Let $\Psi$ be an $\AAa$-localized frame and $w$ an $\AAa$-admissible weight. Then $\left( \Hil_w^\infty , \norm{\Hil_w^\infty}{\cdot} \right)$ is a Banach space, and 
\begin{enumerate} 
\item[(i)] $C_\Psi : \left( \Hil_w^\infty , \norm{\Hil_w^\infty}{\cdot} \right) \rightarrow \ell_w^\infty$ is continuous. 
\item[(ii)] $D_\Psi : \ell_w^\infty \rightarrow \left( \Hil_w^\infty , \norm{\Hil_w^\infty}{\cdot} \right)$ is continuous with \mbox{$\norm{\ell_w^\infty \rightarrow \Hil_w^\infty
}{D_\Psi} \le \norm{\ell_w^\infty \rightarrow \ell_w^\infty}{G_{\tilde
    \Psi,\Psi}}$.} The series $D c = \sum \limits_{k \in K} c_k
\psi_k$ is weak-$^*$ unconditionally convergent. 
\end{enumerate}
\end{lemma} 

\subsection{Duality} \label{sec:dual0}
The associated Banach spaces $\HHmp$ are a generalization of the coorbit spaces in \cite{fegr89} and the modulation spaces \cite{fe06}. We first formulate their duality.
\begin{proposition} \label{sec:dualassocbanach1}  Let $\Psi$ be a $\AAa$-localized frame and $w$ an admissible weight. 
Let $1 \le p < \infty$ and $q$ such that $\frac{1}{p} +\frac{1}{q} = 1$  or $(p,q) = (0,1)$.
Then 
$$ {\left( \HHmp
\right)}' \cong \Hil^q_{1/w}
, $$
where the duality for $f \in \HHmp$ and $h \in \Hil_{1/w}^q$ is given by  
$$\left< f , h \right>_{\HHmp,\HHmpdual} :=  \left< C_{\widetilde \Psi} f, C_{\Psi} h \right>_{\ell_w^p,\ell_{1/w}^q}. $$
\end{proposition} 
\begin{proof} Fix $h \in \Hil_{1/w}^q$. Then, using the duality of $\ell_w^p$ and $\ell_{1/w}^q$, we define a 
sesquilinear form by
 $$ \left< f, h \right>_{\HHmp,\HHmpdual} = 
\left< C_{\widetilde \Psi} f, C_{\Psi} h \right>_{\ell_w^p,\ell_{1/w}^q}, $$ 
 for $f \in \HHmp$.  

Now set $\WWW{h} (f) = \left< f, h \right>_{\HHmp,\HHmpdual}$. 
Then $\left| \WWW{h}(f) \right| \le \norm{\ell_w^p} {C_{\widetilde
    \Psi} f} \norm{\ell_{1/w}^q}{C_{\Psi} h} =
\norm{\ell_{1/w}^q}{C_{\Psi} h} \cdot  \norm{\HHmp}{f}$, since $C_\Psi
h \in \ell_{1/w}^q$ by the norm equivalence
\eqref{sec:syngrammop2}. If $p=0$, we use the estimate
$\left| \WWW{h}(f) \right| \le \norm{\ell_w^\infty} {C_{\widetilde
    \Psi} f} \norm{\ell_{1/w}^1}{C_{\Psi} h} = 
\norm{\ell_{1/w}^1}{C_{\Psi} h} \cdot  \norm{\Hil ^0_w}{f}$, 
Therefore  $\WWW{h} \in \left( \HHmp \right)'$. As a consequence 
$ \WWWfun: 
\HHmpdual \rightarrow 
\HHmpdddd$
is bounded, with $\norm{ \left( \HHmp \right)'}{\WWW{h}} \le
\norm{\ell_{1/w}^q}{C_{\Psi} h}$. 

Conversely, let  $H \in  {\left( \HHmp \right)}'$ and $c \in \llmp$ be
arbitrary with $1\leq p\leq \infty $ or $p=0$. Then $D_\Psi c$ is in $\mathcal{H}^p_w$ and so  $H \left( D_{\Psi} c \right) = \sum \limits_k c_k H(\psi_k)$
with unconditional convergence. 
Therefore the sequence $\left(  H(\psi_k) \right)_k$ is in  $\ell_{1/w}^q$ \cite{conw1}. 
Now define the operator $\VVVfun : \HHmpdddd \rightarrow \HHmpdual  $
by $\VVV{H} = \sum \limits_k \overline{H(\psi_k)} \tilde \psi_k$. 

For $f \in \HHmp$ we have 
$$\WWW{\VVV{ H}} ( f ) = \left< C_{\widetilde \Psi} f , C_\Psi \VVV{ H}\right>_{\ell_w^p, \ell^q_{1/w}} = \sum \limits_l \left< f , \tilde \psi_l\right> 
\overline{\left< \sum \limits_k \overline{H(\psi_k)} \tilde \psi_k  , \psi_l\right>}
= $$
$$ = \sum \limits_k H(\psi_k)  
\left< \sum \limits_l \left< f ,  \tilde \psi_l\right>
 \psi_l , \tilde \psi_k \right>   = 
 \sum \limits_k H(\psi_k)  
\left< f , \tilde \psi_k \right>  =  H \left( \sum \limits_k 
\left< f , \tilde \psi_k \right> \psi_k \right) = H(f).$$
The formal manipulations are justified by the unconditional
convergence of the series involved, by the continuity of $H$, and by
density arguments. Thus 
$ \WWWfun: \HHmpdual \rightarrow 
\HHmpdddd$ is onto.

On the other hand
$$ {\mathcal V} \left({\mathcal W} \left( h \right) \right) =  \sum
\limits_l \overline{\WWW{h}(\psi_l)} \tilde \psi_l =  \sum \limits_l
\overline{\sum \limits_k \left< \psi_l, \tilde \psi_k\right>
  \overline{\left< h, \psi_k\right>} } \tilde \psi_l = $$ 
$$ = \sum \limits_k  \left< h, \psi_k\right> \left( \sum \limits_l \left< \tilde \psi_k , \psi_l \right>  \tilde \psi_l \right)
= \sum \limits_k  \left< h, \psi_k\right> \tilde \psi_k  = h. $$

Therefore $\mathcal{W}$ is invertible. \ep

Similar results appeared in \cite{fegr89} and \cite{dafora05}. 
\begin{remark}
Note that the duality is consistent with the inner product $\left< . ,
  . \right>$ on $\Hil$, see \sloppy  Lemma \ref{sec:dualrel1}. 

Also, note that the isomorphism $ {\big( \HHmp
\big)}' \cong \Hil^q_{1/w}
$ is not an isometric isomorphism.
\end{remark}

By the above result we now have
$\left( \Hil_{1/w}^1 \right)' \cong \left( \Hil_w^\infty,
  \norm{\Hil_w^\infty}{\cdot} \right)$. This yields another proof for
the completeness of 
$\big( \Hil_w^\infty, \norm{\Hil_w^\infty}{\cdot} \big)$. 

\subsubsection{Duality for $\HHmi$}

For $p = \infty$ we can now prove a reconstruction result, as an
extension to Theorem \ref{sec:reconst0}.
\begin{lemma} \label{sec:reconstinf}  Let $\Psi$ be a $\AAa$-localized frame and $w$ an $\AAa$-admissible weight. 
If $f \in \Hil_w^\infty$, then 
 $f = \sum \limits \left< f , \psi_k \right> \tilde \psi_k$ and $f = \sum \limits \left< f , \tilde  \psi_k \right> \psi_k$  with weak-* unconditional convergence in $\sigma\left(\Hil_w^\infty, \Hil ^1_{1/w}\right)$.

Therefore $f = D_{\widetilde \Psi} C_\Psi f = D_{\Psi} C_{\tilde \Psi} f,$ and, 
in particular $D_{\widetilde \Psi}$ and $D_\Psi$ are onto $\HHmi$. 

The norm equivalence \eqref{sec:syngrammop2} is valid for all $1\le p \le \infty$ and $p= 0$:
\begin{equation} \label{sec:syngrammop3}  \tag{\ref{sec:syngrammop2}'}
\frac{1}{\norm{\llmp \rightarrow \llmp}{G_{\widetilde \Psi}}} \norm{\HHmp}{f} \le \norm{\llmp}{C_\Psi f} \le  \norm{\llmp \rightarrow \llmp}{G_{\Psi}} \norm{\HHmp}{f}.
\end{equation}
\end{lemma}
\begin{proof} By above, $\Hil_{1/w}^1$ is the predual of $\Hil_w^\infty$. Now, let $ f \in \Hil_w^\infty$ and $g \in \Hil_{1/w}^1$, then we have 
$$ \Big| \sum \limits_{k \in K} \left< f , \psi_k \right> \left< \tilde \psi_k , g \right> \Big| \le \sum \limits_{k \in K}  \left|  \left< f , \psi_k 
\vphantom{\left|  \left< \tilde \psi_k , g \right> \right|}
 \right> \right|  \left|  \left< \tilde \psi_k , g \right> \right| \le \norm{\ell_w^\infty}{ C_{\Psi} f} \norm{\ell_{1/w}^1}{ C_{\widetilde \Psi}{g}},$$
and the sum of the left-hand side converges absolutely.

 By Lemma \ref{sec:propinft1} $ D_{\tilde \Psi} C_\Psi$ is well-defined on all of $\Hil_w^\infty$. Let $g \in \Hn$, then 
$$ \left< D_{\tilde \Psi} C_\Psi f , g \right> = \lim
\limits_{\tiny \begin{array}{c} H \rightarrow K \\ H
   \,  \mathrm{finite }\end{array}}  \left< \sum \limits_{k \in H} \left< f ,
    \psi_k \right> \cdot \tilde \psi_k , g \right> = \lim
\limits_{\tiny \begin{array}{c} H \rightarrow K \\ H
  \,  \mathrm{finite}\end{array}}   \sum \limits_{k \in H} \left< f , \psi_k
\right> \cdot  \left< \tilde \psi_k , g \right> = \left< f , g
\right>.$$ 
And so $f = D_{\tilde \Psi} C_\Psi f$. 

 The second  reconstruction formula follows by an analogous argument.
The norm equivalence
\eqref{sec:syngrammop3} follows immediately from  the reconstruction formula.
\ep

We can formulate  the compatibility of the  duality relations in the following way.
\begin{lemma} \label{sec:dualrel1}  Let $\Psi$ be a $\AAa$-localized frame and $w$ an admissible weight. For $f \in \HHmp$ and $h \in \HHmidual$ we have 
$$ \left< f , h \right>_{\HHmp, \HHmpdual} = \left< f , h \right>_{\HHmi, \HHmidual}. $$
 
\end{lemma} 
\begin{proof}
The identity follows from the definition of the duality in Proposition \ref{sec:dualassocbanach1}, because
$$\left< f , h \right>_{\HHmp,\HHmpdual} =  \left< C_{\widetilde \Psi}
  f, C_{\Psi} h \right>_{\ell_w^p,\ell_{1/w}^q} = \left< C_{\widetilde
    \Psi} f, C_{\Psi} h \right>_{\ell_w^\infty,\ell_{1/w}^1} = \left< f ,
  h \right>_{\HHmi,\HHmidual}.$$ 
 
\ep 
After clarifying the meaning of the duality brackets, we can now give
the traditional definition of the  coorbit space $\HHmp$ as a  subspace of ``distributions''.
\begin{proposition} \label{sec:distrrep1}
 Let $\Psi$ be an 
 $\AAa$-localized frame and $w$ an admissible weight. 
 For $1 \le p < \infty$ we have
$$\HHmp
= \left\{f \in \HHmi : f = \sum \limits_{k \in K} \left< f , \tilde \psi_k \right>_{\HHmi,\HHmidual} \psi_k 
\mbox{ with }  \left< f , \tilde \psi_k \right>_{\HHmi,\HHmidual} \in \ell_w^p\right\},$$
with unconditional convergence.
\end{proposition}
\begin{proof} 
We combine  the reconstruction formula in Theorem
\ref{sec:origrecon1} with the identities
 $\Big(\left< f , \psi_k \right>_{\HHmi,\HHmidual} \Big)_{k\in K} = \Big(\left< f , \psi_k
\right>_{\HHmp,\HHmpdual} \Big)_{k\in K} \in \ell_w^p$ from
Lemma~\ref{sec:dualrel1} and use the unconditional convergence in
$\HHmp$.  
We obtain 
$$ f = \sum \limits_{k \in K} \left< f , \psi_k
\right>_{\HHmp,\HHmpdual} \tilde \psi_k = \sum \limits_{k \in K}
\left< f , \psi_k \right>_{\HHmi,\HHmidual} \tilde \psi_k .$$ 
\ep 

For $p = \infty$ we can state the following characterization (using
Proposition \ref{sec:dualassocbanach1} and Lemma \ref{sec:reconstinf}): 
\begin{corollary}  Let $\Psi$ be a $\AAa$-localized frame. Let $W$ and $w$  admissible weights satisfying $W \le w$  so that $\HHmi \subseteq \Hil_W^\infty$.
Then for $f \in \Hil_W^\infty$ the following properties are equivalent:
\begin{enumerate}
\item[(i)] $f \in \HHmi$.
\item[(ii)] $\norm{\ell_w^\infty}{ C_{\widetilde \Psi} f } < \infty$.
\item[(iii)] There is a $c \in \ell_w^\infty$, such that $f = \sum
  \limits_{k \in K} c_k \psi_k$ with $C_{\widetilde \Psi} f = G_{\Psi,
    \widetilde \Psi} c$ and 
$$ \norm{\ell_w^\infty}{C_{\widetilde \Psi} f} = \norm{\ell_w^\infty}{G_{\Psi, \widetilde \Psi} c}.$$
\item[(iv)] $f \in \left( \Hil_w^1\right)'.$
\end{enumerate} 
\end{corollary} 

\subsubsection{The chain of Banach spaces $\HHmp$} \label{sec:chainBanachspaces0}

Formulated for Gelfand triples we obtain the following consequence. 
\begin{corollary} \label{sec:gelfand1} Let $\Psi$ be a $\AAa$-localized frame and $w$ an admissible weight with $\inf \limits_{k \in K} w_k > 0$. 
Let $1 \le p < 2$ and $1/p + 1/q = 1$, or $(p,q) = (0,1)$.
Then $\Psi$ is a Gelfand frame for the Gelfand triples 
$$ \Hil_w^p \subseteq \Hil \subseteq \Hil_{1/w}^q,$$
with respect to the duality $\left< C_{\widetilde \Psi} f, C_{\Psi} h \right>_{\ell_w^p,\ell_{1/w}^q}$   
and the sequence spaces $\ell_w^p \subseteq \ell^2 \subseteq \ell_{1/w}^q$. 
\end{corollary} 
\begin{proof}
By Proposition \ref{sec:dualassocbanach1} $ {\left( \HHmp
\right)}' \cong \Hil^q_{1/w}
$. Since $w_k \ge C$ we have, for  $1 \le p \le 2$ and $2 \le q \le \infty$, the following inclusions
$$ \Hil_w^1 \subseteq \Hil_w^p \subseteq \Hil_w^2 \subseteq \Hil \subseteq \Hil_{1/w}^2 \subseteq \Hil_{1/w}^q \subseteq \Hil_{1/w}^\infty,$$
For $q < \infty$, these inclusions are norm-dense, continuous
embeddings (by the corresponding inclusions \eqref{sec:scalseq1} for
sequence spaces), for $q = \infty$, $\Hil$ is w*-dense in
$\Hil_{1/w}^\infty$. Theorem \ref{sec:reconst0} asserts that $\Psi$ is
a Banach frame for $\HHmp$ and $\HHmpdual$. 
\ep
To summarize the ``size'' of the coorbit spaces for $1 \le p \le 2$ and $2 \le q \le \infty$ by Equation \eqref{sec:scalseq1} we obtain the following inclusions: 

\begin{align*} 
\Hn & \subseteq & & \Hil_w^1 & \subseteq & & \Hil_w^p & \subseteq & & \Hil & \subseteq & & \Hil_{1/w}^q & \subseteq & & \Hil_{1/w}^0 & \subseteq & & \Hil_{1/w}^{\infty} \\
 & & & & & &  D_\Psi \bigg\uparrow & & & & & & \bigg\downarrow C_\Psi \\
 \ell^{0 0} & \subseteq & & \ell_w^1 & \subseteq & & \ell_w^p & \subseteq & & \ell^2 & \subseteq & & \ell_{1/w}^q & \subseteq & & \ell_{1/w}^0 & \subseteq & & \ell_{1/w}^{\infty}
\end{align*}

All the inclusions but the last one are in fact norm-dense embeddings;
$\Hn$ is norm dense in $\HHmp$ for $1 \le p < \infty$ and weak*-dense
in $\HHmi$.

Finally we mention that the assumption on the weight serves only  to obtain
a ``small space'' with $p=1$ on the left side of the diagram. By contrast,  if $1/w \subseteq \ell ^2$, then
$\ell ^\infty _w \subset \ell ^2$, and one obtains the Gelfand pair $\Hil^0_w \subseteq \Hil \subseteq \Hil ^1_{1/w}$, which looks a bit
unusual.

\subsubsection{Equivalence result on inclusion of sequence spaces and associated Banach spaces} \label{sec:equivseqassoc0}
Whereas the inclusions of the coorbit spaces $\HHmp$ follow from the
inclusions of the weighted $\ell^p$-spaces, the converse is less obvious and requires more tools.
\begin{theorem}\label{sec:equivclass1}  Let $\Psi$ be an
  $\AAa$-localized norm-bounded frame, i.e.,  $\inf  \limits_k \norm{\Hil}{\psi_k} > 0$. 
 Let $1 \le p_1, p_2 \le \infty$ and let $w_1,w_2$ be admissible weights. Then 
$$ \Hil_{w_1}^{p_1} \subseteq \Hil_{w_2}^{p_2} \Longleftrightarrow \ell_{w_1}^{p_1} \subseteq \ell_{w_2}^{p_2}. $$
\end{theorem} 
\begin{proof} The implication 
$\Longleftarrow$ is clear. \\

Conversely, assume that $\Hil_{w_1}^{p_1} \subseteq \Hil_{w_2}^{p_2}$ and that $\ell_{w_1}^{p_1} \not\subseteq \ell_{w_2}^{p_2}$.

Since $\Psi$ is a frame, by the Feichtinger conjecture
\cite{casazza06} proved in \cite{maspsr15}, $\Psi $ is a finite union of
Riesz sequences.   
In particular, $\Hil$ contains an infinite-dimensional subspace with a Riesz basis $\Psi_0 = \left\{ \psi_k \left| k \in K_0 \right. \right\} \subseteq \Psi$.

The Gram matrix $G_{\Psi_0}$ is  invertible on $\ell^2 (K_0)$ \cite{ole1}.
We extend $G_{\Psi_0}$ to a matrix on $\ell^2(K)$ by defining $G c =
0$ for $c \in \ell^2 ( K \backslash K_0 ) = \ell^2 (K_0)^\bot$. 
Note that $G$ is obtained from the
Gram matrix $G_\Psi$ by setting $(G_\Psi )_{j,k} =0$ for $j,k \not \in
K_0$. Since $\AAa$ is solid, we conclude that  $G$ is in $\AAa$, and since $\AAa$ is closed with respect to the pseudo-inversion \cite{forngroech1}, we also find that $G^\dagger \in \AAa$.
The matrix $G$ possesses the pseudo-inverse $G^\dagger$ with $G^\dagger =
G_{\Psi_0}^{-1} c$ for $c \in \ell^2 (K_0)$ and $G^\dagger c = 0$ for
$c \in \ell^2 (K_0)^\bot$.

Extend  $c \in \ell_{w_1}^{p_1} (K_0) \left\backslash \ell_{w_2}^{p_2} (K_0) \right.$ to a sequence $\tilde c \in \ell_{w_1}^{p_1} (K) \left\backslash \ell_{w_2}^{p_2} (K) \right.$ (by setting $\tilde c_k = 0$ for $k \in K \backslash K_0$) and set $d = G^\dagger \tilde c$.  In particular, $d_k = 0$ for $k \not\in K_0$. Since $G^\dagger \in \AAa$ is bounded on $\ell_{w_1}^{p_1}$, it follows that 
$d \in \ell_{w_1}^{p_1} (K) \backslash \ell_{w_2}^{p_2} (K)$. Therefore $f = \sum \limits_{k \in K} d_k \psi_k \in \Hil_{w_1}^{p_1}$. Furthermore
$$ \left(C_{\Psi_0} f\right)_k = \left< f , \psi_k \right> = \sum \limits_{l \in K_0} d_l \left< \psi_l,\psi_k\right> = \left( G d\right)_k = c_k \text{ for } k \in K_0.$$
If $f \in \Hil_{w_1}^{p_1}$, then $C_{\Psi} f \in \ell_{w_1}^{p_1}.$ 
But ${C_\Psi f}{\big|_{K_0}} = c \not \in \ell_{w_2}^{p_2} (K_0)$, and so ${C_\Psi f} \not \in \ell_{w_2}^{p_2}$. Therefore $f \not \in \Hil_{w_2}^{p_2}$, which is a contradiction.
\ep

Most likely, the statement could be proved without the full strength
of the theorem of Marcus, Spielman, and Srivastava \cite{maspsr15}.

\subsection{Properties of the frame-related operators} \label{sec:framrelop0}

We next summarize the properties of the canonical operators associated to every frame.
We include the statements for $p = \infty$ and discuss the convergence of series expansions. 
As a novelty, we discuss all operators with respect to a frame $\Phi$
in the same localization class, i.e. $\Phi \sim_{\AAa} \Psi$. Being
pedantic, we always consider the operator as a mapping with a
domain. For instance, whereas the synthesis operator is the formula
$D_\Psi c = \sum \limits_{k \in K} c_k \psi_k$ we will use the
notation $D_\Psi^{p,w}$ to denote the synthesis operator on $\HHmp$.

\begin{theorem}\label{sec:locopreint1} Let $\Psi$ be an 
 $\AAa$-localized frame and 
let $w$ be an $\AAa$-admissible weight. Let $1 \le p \le \infty$ and let $1/p + 1/q = 1$ or $(p,q) = (0,1)$. 
\begin{enumerate}
\item[(i)] The analysis operator $C_{\Psi}^{p,w} : \HHmp 
 \rightarrow \ell_w^p$ is given by 
$$C_{\Psi}^{p,w} f = \left( \left< f , \psi_k \right>_{\HHmp, \Hil_{1/w}^q}, k \in K \right).$$
Then 
$C_{\Psi}^{p,w}$ is bounded, one-to-one and has closed range in $\ell_w^p$. 
$C_{\Psi}^{p,w}$ is the restriction of $C_{\Psi}^{\infty,w}$ to $\HHmp$, i.e. $C_{\Psi}^{p,w} = {C_{\Psi}^{\infty,w}}{\big|_{\HHmp}}$.
Furthermore $\range{C_{\Psi}^{p,w}} = \range{C_{\widetilde \Psi}^{p,w}}$ and this  
is a complemented subspace:
\begin{equation} \label{sec:complsubsp0} \ell_w^p = \range{C_{\Psi}^{p,w}} \oplus \kernel{D_{\Psi}^{p,w}}, 
\end{equation}
$\Psi$ is an $\ell ^p _w$-Banach frame for all $\HHmp$ with  bounds
\begin{equation} \label{sec:Banachfram1}
A = \norm{\llmp \rightarrow \llmp}{G_{\widetilde \Psi}^{p,w}}^{-1} \text{ and } B = \norm{\llmp \rightarrow \llmp}{G^{p,w}_{\Psi}}.
\end{equation}

\item[(ii)] The synthesis (or reconstruction) operator $D_{\Psi}^{p,w} : \ell_{w}^p \rightarrow \Hil_{w}^p$ is given by 
$$ D_{\Psi}^{p,w} c = \sum \limits_k c_k \psi_k $$
with unconditional convergence in $\Hil_{w}^p$ for $1 \le p < \infty$
and $p=0$,  and weak*-convergence for $p = \infty$. Then $D_{\Psi}^{p,w}$ is bounded with operator norm $1$ and it maps onto $\HHmp$.
Furthermore $D_{\Psi}^{p,w}$ is the restriction of $D_{\Psi}^{\infty,w}$ to $\HHmp$, i.e. $D_{\Psi}^{p,w} = {D_{\Psi}^{\infty,w}}{\big|_{\HHmp}}$. 
For $p < \infty$
\begin{eqnarray*} 
{\left( D_{\Psi}^{p,w} \right)}^*  & = & C_{\Psi}^{q,1/w} \mbox{, and }\\
{\left( C_{\Psi}^{p,w} \right)}^* & = & D_{\Psi}^{q,1/w}.
\end{eqnarray*}
\item[(iii)] The frame operator $S_\Psi^{p,w}: \Hil_w^p \rightarrow \Hil_w^p$ is defined by
$$ S_\Psi^{p,w} f =  \sum \limits_k \left< f , \psi_k\right>_{\Hil_w^p,\Hil_{1/w}^q} \cdot \, \psi_k =  \sum \limits_k \left< f , \psi_k\right>_{\HHmi,\HHmidual} \cdot \, \psi_k$$ 
with unconditional convergence in $\Hil_w^p$ for $1 \le p < \infty$ and $p = 0$.
It is unconditionally weak*-convergent for $p = \infty$. Furthermore $S_\Psi^{p,w} = D_\Psi^{p,w} C_\Psi^{p,w}$ and $\left( S_\Psi^{p,w} \right)^* = S_\Psi^{q,1/w}$, and $S_{\Psi}^{p,w} = {S_{\Psi}^{\infty,w}}{\big|_{\HHmp}}$. 
The operator $S_\Psi^{p,w}$ is bounded with  bound $\norm{\llmp \rightarrow \llmp}{G^{p,w}_{\widetilde \Psi , \Psi}} \norm{\llmp \rightarrow \llmp}{G^{p,w}_{\Psi}} $. 
It is invertible with inverse $\left(S_\Psi^{p,w}\right)^{-1} = S_{\widetilde \Psi}^{p,w} = {S_{\widetilde \Psi}^{\infty,w}}{\big|_{\HHmp}} = \big({S_{\Psi}^{\infty,w}}{\big|_{\HHmp}}\big)^{-1} $, and is therefore simultaneously invertible on all $\HHmp$.
\item[(iv)] For the Gram matrix $\left(G_{\Psi}\right)_{k,l} = \left< \psi_l, \psi_k \right>_{\Hil}$ (which by the admissibility induces a bounded operator $G_{\Psi}^{p,w} : \ell_{w}^q \rightarrow \ell_w^p$) we have $G_{\Psi}^{p,w} = C_\Psi^{p,w} D_\Psi^{p,w}$ and again 
$G_{\Psi}^{p,w} = {G_{\Psi}^{\infty,w}}{\big|_{\HHmp}}$
The operator $G_{\Psi,\widetilde \Psi}$ is the projection from $\ell_w^p$ on $\range{C_{\Psi}^{p,w}}$.

\end{enumerate}
\end{theorem}
As a consequence of Theorem \ref{sec:locopreint1}, one may now drop
the indices and  write simply and  unambiguously  $C_{\Psi}$, $D_{\Psi}$, $S_{\Psi}$ and $G_{\Psi}$.

We split the proof of Theorem \ref{sec:locopreint1} into shorter lemmata. 
Note that we prove them for an arbitrary frame $\Phi$ that is
localized with respect to the intrinsically localized frame $\Psi$. So
we need the following preparatory result. 
\begin{lemma} \label{sec:crossself1} Let $\Phi$ and $\Psi$ frames with
  $\Phi \sim_{\AAa} \Psi$ and $\Psi \sim_{\AAa} \Psi$ and  let $w$ be an
  $\AAa$-admissible weight.  
Then $\Phi$ is intrinsically localized, and $\widetilde \Phi
\sim_{\AAa} \widetilde \Psi$ and $ \Phi \sim_{\AAa} \widetilde \Psi$.

In particular $\phi_k \in \Hil ^1_{1/w}$ and 
$\left< f, \phi_k \right>_{\HHmp,\HHmpdual} =´\left< f, \phi_k
\right>_{\HHmi,\HHmidual}$ for $f\in \Hil ^p_w$. 
\end{lemma}
\begin{proof}
Since $\widetilde \Psi \sim_{\AAa} \widetilde \Psi$ by
Theorem~\ref{sec:origrecon1}, we may apply  \eqref{sec:nearequiv0} as
follows:
\begin{align*}
   \Phi & \sim_{\AAa}  \Psi, \widetilde \Psi
  \sim_{\AAa} \widetilde \Psi \Rightarrow  \Phi \sim_{\AAa}
  \widetilde \Psi \\
\Phi  & \sim_{\AAa}  \Psi, \widetilde \Psi
  \sim_{\AAa} \Phi \Rightarrow \Phi  \sim_{\AAa}
  \Phi \, .
\end{align*}
As a consequence, the frame $\Phi $ is $\mathcal{A}$-localized and all
results about $\mathcal{A}$-localized frames apply to $\Phi $. 
In particular, Proposition \ref{sec:equivspac0} implies that $\Hil_w^p ( \Psi,
\widetilde \Psi) = \Hil_w^p ( \Phi, \widetilde \Phi )$ with equivalent norms, and we may
write unambiguously $\Hil ^p_w$. 

Furthermore, since the Gram matrix $G_{\widetilde \Psi , \Phi } \in \mathcal{A}$  is
bounded on $\ell ^1_w$, it follows that every row and column of $
G_{\widetilde \Psi , \Phi }$ belongs to $\ell ^1_w$ and likewise to
$\ell ^1_{1/w}$. Consequently, $\phi _k = \sum _{l\in K} \langle \phi
_k , \tilde \psi _l\rangle \psi _l$ is in $\Hil ^1_w \cap \Hil
^1_{1/w}$. Thus the brackets  
$$\left< f, \phi_k \right>_{\HHmp,\HHmpdual} =  \left< f, \phi_k \right>_{\HHmi,\HHmidual}$$
are  well-defined by Lemma~\ref{sec:dualrel1}. 
\ep 
In particular all results shown above are applicable also for $\Phi$,
however,  with  equivalent norms.

\begin{lemma}\label{sec:analoper1}
 Let $\Phi $ and $\Psi$ be  
 $\AAa$-localized frames such that $\Phi \sim_{\AAa} \Psi$.
Let $w$ be an $\AAa$-admissible weight and $1 \le p,q  \le \infty$
with  $1/p + 1/q = 1$ or $(p,q) = (0,1)$. 
Then the  analysis operator $C_{\Phi}^{p,w} : \HHmp 
 \rightarrow \ell_w^p$  given by 
$$C_{\Phi}^{p,w} f = \left( \left< f , \phi_k \right>_{\HHmp, \Hil_{1/w}^q} , k \in K \right)$$
is  bounded, one-to-one and has closed range. 
Furthermore
\begin{equation} \label{sec:syngrammop6}
\frac{1}{ \norm{\llmp \rightarrow \llmp}{G_{\widetilde \Psi, \widetilde \Phi}^{p,w}}} \norm{\HHmp}{f} \le \norm{\llmp}{C_{\Phi}^{p,w} f} \le  \norm{\llmp \rightarrow \llmp}{G^{p,w}_{\Phi, \Psi}} \norm{\HHmp}{f},
\end{equation}
where both sides of the inequality are bounded. The operator $C_{\Phi}^{p,w}$ 
 is the restriction of $C_{\Phi}^{\infty,w}$ to $\HHmp$, i.e. $C_{\Phi}^{p,w} = {C_{\Phi}^{\infty,w}}_{|_{\HHmp}}$.
\end{lemma}

\begin{proof} 
The associated Banach spaces $\HHmp$ coincide for the frames $\Psi$ and $\Phi$.
By Proposition \ref{sec:equivspac0} and Lemma \ref{sec:dualrel1} 
$\left( C_{\Phi}^{p,w} f \right)_k = \left< f , \phi_k \right>_{\HHmp, \HHmpdual}  = \left< f , \phi_k \right>_{\HHmi, \HHmidual}$.
 
\begin{equation*}
\begin{array}{l c c r}
 \norm{\llmp}{C_{\widetilde \Psi}^{p,w} f} &=  \norm{\llmp}{C_{\widetilde \Psi}^{p,w} D_{\widetilde \Phi}^{p,w} C_{\Phi}^{p,w}  f} &\le \norm{\llmp \rightarrow \llmp}{G_{\widetilde \Psi, \widetilde \Phi}^{p,w}} \norm{\llmp}{C_{\Phi}^{p,w}  f} , & \text{ and } \\
 \norm{\llmp}{C_{\Phi}^{p,w} f} &=  \norm{\llmp}{C_{\Phi}^{p,w} D_{\Psi}^{p,w} C_{\widetilde \Psi}^{p,w}  f} &\le \norm{\llmp \rightarrow \llmp}{G_{\Phi, \Psi}^{p,w}} \norm{\llmp}{C_{\widetilde \Psi}^{p,w}  f}.
 \end{array}
 \end{equation*}
By Lemma \ref{sec:crossself1} the Gram matrices $G_{\widetilde \Psi,  \widetilde \Phi}^{p,w}$ and $G_{\Phi, \Psi}^{p,w}$ are in
$\mathcal{A}$ and are therefore bounded on $\ell ^p_w$ for all $p,
1\leq p\leq \infty$. 
By Lemma \ref{sec:dualrel1}  $C_{\Psi}^{p,w} =
{C_{\Psi}^{\infty,w}}_{|_{\HHmp}}$, since  $\Hil ^p_w \subseteq \Hil ^\infty _w$. 
\ep

As $\widetilde \Psi \sim_{\AAa}  \Psi$  the analysis operator
$C_{\widetilde \Psi}^{p,w} : \HHmp 
 \rightarrow \ell_w^p$ is given by  \mbox{$C_{\widetilde \Psi}^{p,w} f
   = \left< f , \tilde \psi_k \right>_{\HHmp, \Hil_{1/w}^q}$}. By
 definition, this particular analysis operator is an isometry. 
 
\begin{lemma}\label{sec:synoper1}
Under the assumptions of Lemma~\ref{sec:analoper1} 
the synthesis operator $D_{\Phi}^{p,w} : \ell_{w}^p \rightarrow
\Hil_{w}^p$ is bounded and onto with operator norm 
$$\norm{\ell_w^p \rightarrow \HHmp}{D_{\Phi}^{p,w}} = \norm{\ell_w^p
  \rightarrow \ell_w^p}{G_{\widetilde \Psi, \Phi}^{p,w}}. $$ 
 It is given by 
$$ D_{\Phi}^{p,w} c = \sum \limits_k c_k \phi_k $$
with unconditional convergence in $\Hil_{w}^p$ for $1 \le p < \infty$
and $p = 0$ and weak*-convergence for $p = \infty$. 
Furthermore, $D_{\Phi}^{p,w}$ is the restriction of $D_{\Phi}^{\infty,w}$ to
$\HHmp$, i.e. $D_{\Phi}^{p,w} = {D_{\Phi}^{\infty,w}}_{|_{\HHmp}}$.  
For $p < \infty$ we have  

\begin{equation} \label{sec:dualoper1} 
{\left( D_{\Phi}^{p,w} \right)}^*  =  C_{\Phi}^{q,1/w} \quad  \text{
  and }  \quad 
{\left( C_{\Phi}^{p,w} \right)}^* = D_{\Phi}^{q,1/w}.
\end{equation}
\end{lemma}
\begin{proof} 
By Lemma \ref{sec:crossself1} and Lemma \ref{sec:reconstinf}
$D^{p,w}_\Phi$ is bounded and onto $\HHmp (\Phi, \widetilde \Phi ) =
\HHmp (\Psi , \widetilde \Psi )$ (see also  Proposition
\ref{sec:equivspac0}).  Since $D^{p,w}_\Phi = G_{\Phi,\widetilde \Psi}
D^{p,w}_\Psi$, $D^{p,w}_\Phi $  is bounded on $\ell ^p_w$.
The unconditional convergence of $D_{\Phi}^{p,w} c = \sum \limits_k
c_k \phi_k $. is shown  as in Lemma \ref{sec:propinft1}. 
Since $\ell ^p_w \subset \ell ^\infty _w$, it is clear that 
$D_{\Phi}^{p,w} = {D_{\Phi}^{\infty,w}}_{|_{\HHmp}}$.

For the adjoint operator  let $c \in \ell_w^p$, and $f \in
\Hil_{1/w}^q \simeq (\Hil ^p_w)^\prime $. Then 
$$\left< D_{\Phi}^{p,w} c , f\right>_{\Hil_{w}^p,\Hil_{1/w}^q} =
\big\langle  \sum \limits_k c_k  \phi_k , f\big\rangle _{\Hil_{w}^p,\Hil_{1/w}^q}
=$$  
$$ =\sum \limits_k c_k \left< \phi_k , f\right>_{\Hil_{w}^p,\Hil_{1/w}^q} = \left< c , C_{\Phi}^{q,1/w} f\right>_{\ell_{w}^p,\ell_{1/w}^q}, $$
where the change of order is justified because $c \in \ell_w^p$, $f \in \HHmpdual$ and by Lemma \ref{sec:analoper1}. 

The operator norm is
$$ \norm{\ell_w^p \rightarrow \HHmp}{D_\Phi^{p,w}} = \sup \limits_{\norm{\ell_w^p}{c} = 1}  \norm{\HHmp}{D_\Phi^{p,w} c}  = \sup \limits_{\norm{\ell_w^p}{c} = 1} \norm{\ell_w^p}{G_{\widetilde \Psi, \Phi}^{p,w} c} = \norm{\ell_w^p \rightarrow \ell_w^p}{G_{\widetilde \Psi, \Phi}^{p,w}c}.$$

\ep

Clearly, by above, $ \norm{\ell_{w}^p \rightarrow \HHmp}{D_\Phi^{p,w}}
= \norm{\HHmpdddd \rightarrow \ell_{1/w}^q}{C_\Phi^{q,1/w}}$.  In
general we have  $\norm{\HHmpdddd \rightarrow
  \ell_{1/w}^q}{C_\Phi^{q,1/w}}\neq \norm{\HHmpdual
  \rightarrow \ell_{1/w}^q}{C_\Phi^{q,1/w}}$, because the isomorphism
between  $(\Hil ^p_w)^\prime$ and $\Hil ^q_{1/w}$ of Proposition
\ref{sec:dualassocbanach1}  is not an isometry.

\begin{lemma}\label{sec:gramoper1} 
Let  $\Psi$, $\Phi$ and $\Xi$ be frames with $\Psi \sim_{\AAa} \Psi$,
$\Phi \sim_{\AAa} \Psi$ and $\Xi \sim_{\AAa} \Psi$, and 
let $w$ be an $\AAa$-admissible weight. Let $1 \le p \le \infty$ and let $1/p + 1/q = 1$ or $(p,q) = (0,1)$.  
The cross-Gram matrix $G_{\Phi,\Xi}$ with entries $\left(G_{\Phi,\Xi}\right)_{k,l} = \left< \xi_l, \phi_k \right>_{\Hil}$ induces a bounded operator $G_{\Phi,\Xi}^{p,w} : \ell_{w}^p \rightarrow \ell_w^p$ and factors as
$$G_{\Phi,\Xi}^{p,w} = C_\Phi^{p,w} D_\Xi^{p,w}, $$
and so $\left( G_{\Phi,\Xi}^{p,w}\right)^* =
G_{\Xi,\Phi}^{q,1/w}$. The Gram matrix $G_{\Phi,\Xi}^{p,w}$ is the
restriction ${G_{\Phi,\Xi}^{\infty,w}}\Big| _{{\HHmp}}$. 

Furthermore  $\range{G_{\Phi,\Xi}^{p,w}}= \range{C_\Phi^{p,w}}$ and $\kernel{G_{\Phi,\Xi}^{p,w}} = \kernel{D_\Xi^{p,w}}$. 
The Gram matrix $G_{\Phi,\Xi}^{p,w}$ is a bijective mapping from $\range{C_{\widetilde \Xi}^{p,w}}$ onto $\range{C_\Phi^{p,w}}$. 
 
For $\Xi = \widetilde \Phi$ the Gram matrix satisfies 
\begin{equation} \label{sec:gramselfadj1}
G_{\Phi,\widetilde \Phi} = G_{\widetilde \Phi,\Phi} = G_{\Phi,\widetilde \Phi}^*, 
\end{equation} 
and $G_{\Phi,\widetilde \Phi}$ is a bounded  projection from
$\ell_w^p$ on the range of $C_\Phi^{p,w}$ with kernel $\kernel{D_{\Phi}^{p,w}}$. 
In addition, $\range{C_\Phi^{p,w}} = \range{C_{\widetilde \Phi}^{p,w}}$ and
\begin{equation*} \label{sec:complsubsp1} 
\ell_w^p = \range{C_{\Phi}^{p,w}} \oplus \kernel{D_{\Phi}^{p,w}}, 
\end{equation*}

Therefore,  we have
\begin{equation} \label{sec:normsynthgramop2}
\norm{\ell_w^p \rightarrow \HHmp}{D_{\Phi}^{p,w}} = \norm{\ell_w^p
  \rightarrow \ell_w^p}{G_{\widetilde \Phi, \Phi}} \ge 1. 
\end{equation} 
\end{lemma}
\begin{proof} 

For $c \in \ell^{00}$ we have 
$$ \left( G_{\Phi,\Xi}^{p,w} \cdot c \right)_l = \sum \limits_{k} \left(G_{\Phi,\Xi}\right)_{l,k} c_k =  \sum \limits_{k} \left< \xi_k , \phi_l\right>   c_k =  \left( C_\Phi D_\Xi c \right)_l . $$ 
Therefore $G_{\Phi,\Xi}^{p,w} = C_\Phi^{p,w} D_\Xi^{p,w}$ on
$\ell^{00}$. Since both sides are bounded operators on $\ell ^p_w$
($G_{\Phi,\Xi}^{p,w}$ because $\Phi \sim _{\mathcal{A}} \Xi$), the
factorization can be extended from the dense subspace $\ell ^{00}$ to
$\ell ^p_w$ for $p<\infty $.

By Lemma \ref{sec:analoper1}  $C_\Phi^{p,w}$ is one-to-one, and
therefore $\kernel{G_{\Phi,\Xi}^{p,w}} =
\kernel{D_\Xi^{p,w}}$. Likewise, by 
Lemma \ref{sec:synoper1} $D_\Xi^{p,w}$ is onto $\Hil ^p_w$, and therefore
$\range{G_{\Phi,\Xi}^{p,w}} = \range{C_\Phi^{p,w}}$.  
Since  $G_{\Phi,\Xi}^{p,w} C_{\widetilde \Xi}^{p,w} f = C_{\Phi}^{p,w}
f$, the Gram matrix $G_{\Phi,\Xi}^{p,w}$  induces a bijective mapping from
$\range{C_\Xi^{p,w}}$ onto $\range{C_{\Phi}^{p,w}}$. (Compare to the
'frame transformation' in \cite{xxlfinfram1}.)

If $\Xi = \widetilde \Phi$, then 
$$ (G_{\Phi,\widetilde \Phi})_{k,l} = \left< \phi_l, 
  \tilde \phi_k \right>_{\Hil} = \left< \phi_l, S^{-1} \phi_k
\right>_{\Hil} = \left< S^{-1} \phi_l, \phi_k \right>_{\Hil} =
(G_{\widetilde \Phi, \Phi} )_{k,l}, $$
and for the entries of the adjoint matrix
$$ (G_{\Phi,\widetilde \Phi}^*)_{k,l} =
\overline{(G_{\Phi,\widetilde \Phi})}_{l,k} =
\overline{\left< \phi_k, S^{-1} \phi_l \right>_{\Hil}} = \left< S^{-1}
  \phi_l, \phi_k \right>_{\Hil}= (G_{\Phi,\widetilde
    \Phi})_{k,l}  ,$$
and \eqref{sec:gramselfadj1} is verified. 

Since $D^{p,w}_{\tilde \Phi  } C^{p,w}_\Phi = \mathrm{Id}_{\Hil ^p_w}$
by \eqref{sec:reconform1}, we obtain
$$
(G_{\Phi,\widetilde \Phi}^{p,w})^2 = C^{p,w}_\Phi D^{p,w}_{\tilde \Phi } C^{p,w}_\Phi D^{p,w}_{\tilde
  \Phi } = C^{p,w}_\Phi D^{p,w}_{\tilde \Phi } = G_{\Phi,\widetilde
  \Phi}^{p,w} \, .
$$
Thus $G_{\Phi,\widetilde \Phi}^{p,w} $ is a projection operator on
$\ell ^p_w$ with range $\range{C_\Phi^{p,w}}$ and kernel
$\kernel{D_{\tilde \Phi } ^{p,w}}= \kernel{D_{\Phi } ^{p,w}} $.
In particular, $\| G_{\Phi,\widetilde \Phi}^{p,w} \| \geq 1$ and
\eqref{sec:normsynthgramop2} follows. 

  Since
$G_{\Phi,\widetilde \Phi}^{p,w} = G_{\widetilde \Phi, \Phi}^{p,w} $ we get $\range{G_{\Phi,\widetilde \Phi}^{p,w}} = \range{G_{\widetilde \Phi, \Phi}^{p,w}}$ and by above 
we also have $\range{C_\Psi^{p,w}} = \range{C_{\widetilde \Psi}^{p,w}}$.
By the projection property, using \cite[Theorem III.13.2]{conw1},
$\range{C_{\Phi}^{p,w}}$ and $\kernel{D_{\Phi}^{p,w}}$ are therefore
complementary subspaces. 
\ep

By combining all properties of $C_\Psi ^{p,w}$ and $D_\Psi ^{p,w}$,
we finally obtain the following list of properties for the frame
operator $S_\Phi ^{p,w} = D_\Phi ^{p,w} C_\Phi ^{p,w} $.
\begin{lemma}\label{sec:framoper1} 
The  frame operator $S_\Phi^{p,w}: \Hil_w^p \rightarrow \Hil_w^p$ is defined by
$$ S_\Phi^{p,w} f =  \sum \limits_k \left< f ,
  \phi_k\right>_{\Hil_w^p,\Hil_{1/w}^q} \cdot \, \phi_k =  \sum
\limits_k \left< f , \phi_k\right>_{\HHmi,\HHmidual} \cdot \, \phi_k   $$ 
with unconditional convergence in $\Hil_w^p$ for $p < \infty$ and
weak*-unconditional convergence for $p = \infty$. 
The frame operator satisfies the identities
$S_\Phi^{p,w} = D_\Phi^{p,w} C_\Phi^{p,w}$,
$\left(S_\Phi^{p,w}\right)^* = S_\Phi^{q,1/w}$, and $S_{\Phi}^{p,w} =
{S_{\Phi}^{\infty,w}}_{|_{\HHmp}}$,  and is  bounded on all $\Hil
^{p,w}$ with operator norm
\begin{equation*}
\begin{array}{l c l c c r}
\norm{\HHmp \rightarrow \HHmp}{S_{\Phi}^{p,w}}  & \le & \norm{\ell_{w}^p \rightarrow \ell_{w}^p}{G_{\widetilde \Psi,\Phi}^{p,w}} & \cdot & \norm{\ell_{w}^p \rightarrow \ell_{w}^p}{G_{\Phi,\Psi}^{p,w}}.&
\end{array}
\end{equation*}
Furthermore, $S_\Phi^{p,w}$ is simultaneously invertible on all
$\HHmp$  with inverse
$\left(S_\Phi^{p,w}\right)^{-1} = S_{\widetilde \Phi}^{p,w}$. 
\end{lemma}

\section{Galerkin Matrix representation of operators with localized frames} \label{sec:matreplocalfr0}

For a numerical treatment of operator equations  one often uses redundant frame representations for the Galerkin discretization. 
Such discretizations have been formulated for wavelet frames in
\cite{Stevenson03} and for Gabor frames in \cite{GroechSjoe1}. The
formalism for general (Hilbert space) frames has been introduced in \cite{xxlframoper1}. 

For  localized frames  we  formally define the relation between
operators and matrices as follows.  
\begin{definition} Let $\Psi$, $\Phi$ and $\Xi$ be frames with $\Psi
  \sim_{\AAa} \Psi$, $\Phi \sim_{\AAa} \Psi$ and $\Xi \sim_{\AAa}
  \Psi$.  Let $w_1, w_2$ be $\AAa$-admissible weights and  $1 \le
  p_1,p_2 \le \infty$ or $p_1,p_2 = 0$. Let $q_1,q_2$ be the dual
  indices defined as usual.  
\begin{enumerate} \item[(i)] For the bounded linear operator $O :
  \Hil_{w_1}^{p_1} \rightarrow \Hil_{w_2}^{p_2}$ define the matrix $
  {\mathcal M}_{(\Phi , \Xi)}$ 
by 
$$ {\left( {\mathcal M}_{(\Phi , \Xi)} \left( O \right) \right)}_{k,l} = 
\left<O \xi_l, \phi_k \right>_{\Hil^{p_2}_{w_2},\Hil_{1/w_2}^{q_2}} =
\left<O \xi_l, \phi_k \right>_{\Hil^{\infty }_{w_2},\Hil_{1/w_2}^{1}} 
. $$
We call ${\mathcal M}_{(\Phi , \Xi)} \left( O \right)$ the \emph{(Galerkin) matrix of $O$ with respect to $\Phi$ and $\Xi$}.
\item[(ii)] For the matrix $M$ that induces a bounded operator in
  $\BL{\ell_{w_1}^{p_1},\ell_{w_2}^{p_2}}$ define 
$ {\mathcal O}_{(\Phi , \Xi)}
: \BL{\ell_{w_1}^{p_1},\ell_{w_2}^{p_2}} \rightarrow  \BL{\Hil_{w_1}^{p_1},\Hil_{w_2}^{p_2}}$ by
\begin{equation} \label{sec:opematr2} \left( \mathcal{O}_{(\Phi ,
      \Xi)} \left( M \right) \right) h = \sum \limits_k  \Big( \sum
    \limits_j M_{k,j} \left<h, \xi_j\right> \Big) \phi_k,  
\end{equation}
for $h \in \Hil_{w_1}^{p_1}$.
We call ${\mathcal O}_{(\Phi , \Xi)} \left( M \right)$ the \emph{operator of $M$ with respect to $\Phi$ and $\Xi$}. 
\end{enumerate} 
\end{definition} 

\begin{theorem} \label{sec:matbyfram1} Assume that $\Phi, \Psi$ and
  $\Xi $ are $\mathcal{A}$-localized frames in $\Hil$  satisfying $\Phi
  \sim_{\AAa} \Psi$ and $\Xi \sim_{\AAa} \Psi$.  Let $w_1, w_2$ be
  $\AAa$-admissible weights, let $1 \le p_1,p_2 \le \infty$ or 
  $p_1,p_2 = 0$  with dual
  indices  $q_1,q_2$. 
\begin{enumerate} \item[(i)] 
If  $O \in \BL{\Hil_{w_1}^{p_1},\Hil_{w_2}^{p_2}}$, then $\mathcal{M}_{(\Phi , \Xi)}  ( O ) \in
\BL{\ell_{w_1}^{p_1},\ell _{w_2}^{p_2}}$,   and  we have 
\begin{equation}
\label{sec:matrixrepA}
\norm{\ell _{w_1}^{p_1} \rightarrow \ell _{w_2}^{p_2}}{{\mathcal M}_{(\Phi , \Xi)}
 \left( O \right)} \le 
\norm{\ell_{w_2}^{p_2} \rightarrow \ell_{w_2}^{p_2}}{G
_{\Phi, \Psi}} 
 \norm{\ell_{w_1}^{p_1} \rightarrow \ell_{w_1}^{p_1}}{G_{\widetilde \Psi, \Xi}
} 
\norm{\Hil_{w_1}^{p_1} \rightarrow \Hil_{w_2}^{p_2}}{O}.
\end{equation}
Furthermore, 
$$
\mathcal{M}_{(\Phi , \Xi)} ( O )  =   C_{\Phi} \circ O
\circ D_{\Xi}  \, . 
$$
\item[(ii)] 
If $M\in \BL{\ell_{w_1}^{p_1},\ell _{w_2}^{p_2}}$,   then $\mathcal{O}_{(\Phi ,
  \Xi)}(M) \in  \BL{\Hil_{w_1}^{p_1},\Hil _{w_2}^{p_2}}$, and 
$$ \mathcal{O}_{(\Phi , \Xi)} (M) = D_{\Phi} \circ M \circ C_{\Xi},  $$
and
\begin{equation}\label{sec:matrixrepB}
\norm{\Hil_{w_2}^{p_2} \rightarrow \Hil_{w_1}^{p_1}}{\mathcal{O}_{(\Phi , \Xi)}\left( M \right)} \le  
 \norm{\ell_{w_1}^{p_1} \rightarrow \ell_{w_1}^{p_1}}{G_{\widetilde \Psi, \Phi}} 
 \norm{\ell_{w_2}^{p_2} \rightarrow \ell_{w_2}^{p_2}}{G_{\Xi, \Psi}} 
\norm{\ell_{w_2}^{p_2} \rightarrow \ell_{w_1}^{p_1}}{M}.
\end{equation}
\end{enumerate} 
\end{theorem}
\begin{proof}
This result follows directly from the results in Section
\ref{sec:localfr0}.  For example, let  $c = (c_k) \in
\ell_{w_1}^{p_1}$, then 
$$ \left({\mathcal M}_{(\Phi , \Xi)} \left( O \right) c\right)_l  =  \left( C
_{\Phi} \circ O \circ D
_{\Xi} c \right)_l = $$
$$ = \left(C_{\Phi} \left( \sum \limits_{k \in K} c_k O \xi_k \right)\right)_l = \sum \limits_{k \in K} c_k \left< O \xi_k , \phi_l \right>_{\Hil_{w_2}^{p_2},\Hil_{1/w_2}^{q_2}}.$$
\ep 

Using tensor products and the results in Section \ref{sec:localfr0},
it is  easy to extend the results in \cite{xxlframoper1}. 
For   $g \in X_1'$, $f \in X_2$  the \emph{tensor product} $f \otimes
g$ is defined as the rank-one operator  from $X_1$ to $X_2$ by
$ \left( f \otimes g \right) (h) = \left< h, g \right>_{X_1,X_1 '}
f$. We will use
$\left( f \otimes g \right) (h) = \left< h, g \right>_{\HHmp ,
  \HHmpdual} f$ for $h\in  \HHmp$. 

\begin{proposition} \label{sec:propmatropfram1} Let $\Psi$ and $\Phi$
  be $\AAa$-localized frames in $\Hil$ satisfying $\Phi \sim
  _{\mathcal{A}} \Psi $ and
  $w_1, w_2$ be $\AAa$-admissible weights and $1 \le p_1,p_2 \le
  \infty$ or $p_1,p_2= 0$. Then the factorization 
$$ \left( \mathcal{O}_{(\Phi , \Psi)}\circ \mathcal{M}_{(\widetilde{\Phi}, \widetilde{\Psi})}\right)  = \identity{} = \left( {\mathcal O_{(\widetilde{\Phi}, \widetilde{\Psi})}
\circ \mathcal{M}_{(\Phi , \Psi)}
}\right) ,$$
holds for every space of bounded operators $\BL{\Hil_{w_1}^{p_1},\Hil_{w_2}^{p_2}}$.

Therefore every $O \in \BL{\Hil_{w_1}^{p_1},\Hil_{w_2}^{p_2}}$ 
possesses the representation 
\begin{equation} \label{sec:operrep2} O = \sum \limits_{k,j} \left<O \tilde{\psi}_j, \tilde{\phi}_k \right>  \phi_k \otimes {\psi}_j = \sum \limits_{k,j} \left<O {\psi}_j, {\phi}_k \right>  {\tilde \phi}_k \otimes {\tilde \psi}_j ,
\end{equation}
and both expansions  converge unconditionally in the strong operator
topology  (respectively weak-* unconditionally if either  $p_1 = \infty$ or
$p_2 = \infty$). 
\end{proposition}
\begin{proof} 
By Theorem \ref{sec:matbyfram1} for an $O $ in any $\BL{\Hil_{w_1}^{p_1},\Hil_{w_2}^{p_2}}$ we have
\begin{equation}\label{sec:operatid1}
\left( \mathcal{O}_{(\Phi , \Psi)}\circ \mathcal{M}_{(\widetilde{\Phi}, \widetilde{\Psi})}\right) O = D_{\Phi} \left( C_{\widetilde \Phi} \,O \,  D_{\widetilde \Psi} \right) C_{\Psi}  = O, 
\end{equation}
using the reconstruction formulas in Theorem \ref{sec:origrecon1} and Lemma  \ref{sec:reconstinf}. 

The representation in \eqref{sec:operrep2} converges in the strong operator topology by Theorem \ref{sec:matbyfram1}.
\ep

As in the Hilbert space setting \cite{xxlframoper1} we get the
following decomposition. 
\begin{proposition}
Let $\Psi, \Phi$ and $\Xi$ be $\AAa$-localized frames in $\Hil$
satisfying $\Phi \sim _{\mathcal{A}} \Psi , \Xi \sim _{\mathcal{A}}
\Psi $.    Let  $w_1, w_2, w_3$ be $\AAa$-admissible weights and 
$1 \le p_1,p_2 \le
  \infty$ or $p_1,p_2= 0$.
Then 
 for $O_1 : \Hil_{w_1}^{p_1} \rightarrow \Hil_{w_2}^{p_2}$ and $O_2 : \Hil_{w_3}^{p_3} \rightarrow \Hil_{w_1}^{p_1}$, we have 
$${\mathcal M}_{(\Phi , \Psi)} \left( O_1 \circ O_2  \right) = {\mathcal M}_{(\Phi , \Xi)} \left( O_1 \right) \circ {\mathcal M}_{(\widetilde \Xi , \Psi)} \left( O_2 \right).$$
\end{proposition}
\begin{proof}
The statement follows from the factorization 
$${\mathcal M}_{(\Phi , \Psi)} \left( O_1 \circ O_2  \right) = C_\Phi \, O_1 \, O_2 \, D_\Psi = $$
$$ = C_\Phi \, O_1 \, D_\Xi \, C_{\widetilde \Xi} \, O_2 \, D_\Psi = {\mathcal M}_{(\Phi , \Xi)} \left( O_1 \right) \circ {\mathcal M}_{(\widetilde \Xi , \Psi)} \left( O_2 \right).$$
\ep 

Then we get an extension of results in \cite{xxlriek11,xxlrieck11} to
coorbit spaces. 
\begin{lemma}\label{sec:boundedopequiv1} Let $\Psi$ ,and $\Phi$  be
  $\AAa$-localized sequence in $\Hil$ satisfying $\Phi \sim
  _{\mathcal{A}} \Psi $, $w_1, w_2$ be $\AAa$-admissible weights and $1
    \le p_1,p_2 \le \infty$ or $p_1,p_2 = 0$. Let $O$ be a linear
    operator from   $\Hil ^{00}$ into $\Hil
    ^\infty _{w_2}$.     Then 
$$ O \in \BL{\Hil_{w_1}^{p_1},\Hil_{w_2}^{p_2}} \Longleftrightarrow
{\mathcal M_{(\Psi,{\Phi})} (O)} \in
\BL{\ell_{w_1}^{p_1},\ell_{w_2}^{p_2}},$$ 
\end{lemma}
\begin{proof} The implication $\Rightarrow$  is stated  in Theorem \ref{sec:matbyfram1}(i).

For the converse, let $O$ be a linear 
operator from  $\Hil ^{00} $  to $\Hil_{w_2}^{p_2}$ such that 
$\mathcal{M}_{(\Psi , \Phi )}(O)= C^{p_2,w_2}_{\Phi} \circ O \circ D^{p_1,w_1}_{\Psi}$ is bounded. 
Then $\range{D_\Psi} \subseteq \domain{O}$ and therefore $O$ is defined everywhere.
Since 
$$ O =  D^{p_2,w_2}_{\widetilde \Phi} \circ C^{p_2,w_2}_{\Phi} \circ O \circ D^{p_1,w_1}_{\Psi} \circ C^{p_1,w_1}_{\widetilde \Psi}, $$
the operator $O$ is also bounded.
\ep

\subsection{Characterization of operator classes}\label{sec:charaop0}

Combining the matrix representation of an operator  with the
well-known characterizations of boundedness of operators between
$\ell^p$-spaces \cite{Maddox:101881} we obtain 
criteria for the boundedness of operators between certain coorbit
spaces. 

For the description we recall the following norms on  infinite
matrices. Considering an index set $K = L \times N$, we can define
discrete mixed norm spaces \cite{gr01}, i.e.,  
$$\ell^{p_1,p_2}_{w} = \left\{ M = \left(M_{l,n}\right)_{l \in L, n \in N} \left| \norm{\ell^{p_1,p_2}_{w}}{M} := \left( \sum \limits_l \left( w_{l,n} \left| M_{l,n} \right|^{p_2}  \right)^{p_1/p_2} \right)^{1/p_2} < \infty\right. \right\}.$$
In particular we consider weights $w = w^{(1)} \otimes w^{(2)}$ with $w_{k,l} = w^{(1)}_k \cdot w^{(2)}_l$. 

\begin{proposition} \label{prrr} Let $\Psi$ and $\Phi$  be $\AAa$-localized
  sequence in $\Hil$ satisfying $\Phi \sim _{\mathcal{A}} \Psi $,
  $w^{(i)}$ be $\AAa$-admissible weights and $1 \le p_i < \infty$ for
  $i = 1,2$.   
Let $O$ be a linear operator from $\Hil ^{00}$ to  $\Hil_{w_2}^\infty$, and $M = {\mathcal M}_{(\Phi , \Xi)} \left( O \right)$. Then 

\begin{equation} \label{c:2}
\begin{array}{l c l}  O \in \BL{\Hil_{w^{(1)}}^\infty,
    \Hil_{w^{(2)}}^{\infty}} & \Longleftrightarrow &  O \in
  \BL{\Hil_{w^{(1)}}^0, \Hil_{w^{(2)}}^{\infty}}\\ 
& \Longleftrightarrow & \sup \limits_{k} \sum \limits_l \left|
  w^{(2)}_k \cdot \frac{1}{w^{(1)}_l}  \left<O \psi_l, \phi_k \right>
\right| < \infty, \\ 
 & \Longleftrightarrow &   
M \in \ell_{1/w^{(2)}  \otimes w^{(1)} }^{\infty,1} .
\end{array}
\end{equation}

\begin{equation} 
\begin{array}{l c l}  O \in \BL{\Hil_{w^{(1)}}^\infty,
    \Hil_{w^{(2)}}^0} & \Longleftrightarrow & \lim \limits_{k} \sum
  \limits_l \left| w^{(2)}_k \cdot \frac{1}{w^{(1)}_l}  \left<O
      \psi_l, \phi_k \right> \right| = 0. 
\end{array}
\end{equation}

\begin{equation}
\begin{array}{l c l} 
O \in \BL{\Hil_{w^{(1)}}^{1},\Hil_{w^{(2)}}^{\infty}} &
\Longleftrightarrow & \sup \limits_{k,l} \left| w^{(2)}_k \cdot
  \left<O \psi_l, \phi_k \right> \cdot \frac{1}{w^{(1)}_l}  \right| <
\infty  ,\\
  & \Longleftrightarrow &   
M \in \ell_{1/w^{(2)}  \otimes w^{(1)} }^{\infty,\infty} .
\end{array}
\end{equation}

\begin{equation}
\begin{array}{l c l} 
O \in \BL{\Hil_{w^{(1)}}^1, \Hil_{w^{(2)}}^{p}} & \Longleftrightarrow &
  \sup \limits_{l} \sum \limits_k \left| w^{(2)}_k \cdot
    \frac{1}{w^{(1)}_l}  \left<O \psi_l, \phi_k \right> \right|^p  <
  \infty, \\ 
& \Longleftrightarrow &  M^* \in \ell_{1/w^{(2)} \otimes w^{(1)}}^{p,\infty}\\
\end{array}
\end{equation}
\begin{equation}
\begin{array}{l c l} 
O \in \BL{\Hil_{w^{(1)}}^\infty, \Hil_{w^{(2)}}^{1}} &
\Longleftrightarrow  & \sup \limits_{E \text{ finite }}  \sum
\limits_l \left| \sum \limits_{k \in E} w^{(2)}_k \cdot
  \frac{1}{w^{(1)}_l}  \left<O \psi_l, \phi_k \right> \right| < \infty
. 
\end{array}
\end{equation}

\begin{equation}
\begin{array}{l c l} 
O \in \BL{\Hil_{w^{(1)}}^2, \Hil_{w^{(2)}}^{2}} & \Longleftrightarrow  & \left\{
\begin{array}{c} \text{For } \breve{M}_{k,l} = w^{(2)}_k \cdot \frac{1}{w^{(1)}_l}  \left<O \psi_l, \phi_k \right> \text{ we have:}  \\
\sum \limits_l \left| \sum \limits_{k \in E} \breve{M}_{k,l} \right|^2 < \infty , \\ \left(\breve{M}^* \breve{M}\right)^n \text{is defined for all } n=1,2, \dots \\ \sup \limits_n \sup \limits_i \left[ \left( \left(\breve{M}^* \breve{M}\right)^n \right)_{i,i} \right]^{1/n} = K < \infty. 
\end{array} \right.
\end{array}
\end{equation}

\end{proposition}
\begin{proof}
The conditions on the matrix $M$ are variations of the well-known
Schur test. For instance, $M$ is bounded from $\ell ^\infty $ to $\ell
^\infty $, if and only if the row sums are uniformly bounded, i.e.,
$\sup _k \sum _l |M_{k,l}| <\infty $. A convenient reference for
Schur's test is \cite{Maddox:101881}. 

To also include weights, we proceed as follows. 
Let $D_j c = (w_k^{(j)} c_k)_{k\in K}$ be the multiplication operator
with weight $w^{(j)}, j=1,2$. Then $D_j$ is an isometric isomorphism
from $\ell ^p_{w^{(j)}}$ onto $\ell ^p$. 

Therefore a matrix $M$ is bounded from $\ell_{w^{(1)}}^{p_1}$ into
$\ell_{w^{(2)}}^{p_2}$, if and only if $\breve{M}= D_2 M D_1^{-1}$  is bounded from
$\ell^{p_1}$ into $\ell^{p_2}$.  
We now apply the  boundedness characterizations  in
\cite{Maddox:101881} to $ D_2 \mathcal{M}_{(\Psi  , \Phi )}(O)
D_1^{-1} $. For example,  Lemma~\ref{sec:boundedopequiv1} says that 
$ O \in \BL{\Hil_{w_1}^{\infty},\Hil_{w_2}^{\infty}} \Longleftrightarrow
{\mathcal M_{(\Psi,{\Phi})} (O)} \in
\BL{\ell_{w_1}^{\infty },\ell_{w_2}^{\infty}}$, which in turn is equivalent
to saying that $ \breve{M}= D_2 \mathcal{M}_{(\Psi,{\Phi})} (O) D_1^{-1} \in \BL{\ell ^\infty, \ell
  ^\infty }$. Since $\breve{M}_{k,l} = w^{(2)}_k \cdot
    \frac{1}{w^{(1)}_l}  \left<O \psi_l, \phi_k \right>$, condition
    \eqref{c:2} follows from cite[Theorem 2.6]{Maddox:101881}. The other
    characterization follow in the same way from \cite[Theorem 2.12,2.13(a),2.13(b),2.14]{Maddox:101881} and \cite{cron71}, respectively.  
\ep

In concrete applications one uses only the sufficient conditions for
boundedness and checks that the matrix $\mathcal{M}_{(\Psi,{\Phi})}
(O)$ satisfies the conditions of Schur's test.  
\begin{corollary} Let $\Psi$ and $\Phi$  be $\AAa$-localized frame in
  $\Hil$ satisfying $\Phi \sim _{\mathcal{A}} \Psi $, 
  and $1 \le p \leq  \infty$.  Let $O$ be a
  linear operator from $\Hil ^{00}$ to  $\Hil_w^\infty$.   
If  
$\sup \limits_{k} \sum \limits_l \left| \left<O \psi_l, \phi_k \right> \right| < \infty,$ and $\sup \limits_{l} \sum \limits_k \left| \left<O \psi_l, \phi_k \right> \right| < \infty,$ 
then $O \in \BL{\Hil ^p, \Hil ^p}$.
\end{corollary}

For an example of how these abstract results are applied in analysis,
we refer to the investigation of the boundedness of pseudodifferential
operators with the help of Gabor frames in~\cite{GroechSjoe1,grrz08}.\\

For further reference, we remark that the results in
Section~\ref{sec:matreplocalfr0} do not use the full power of
intrinsic localization, but remain true under weaker assumptions. In
fact, we have only used the norm 
equivalences for the analysis operators of two frames $\Phi $ and
$\Psi $ and their duals $\tilde \Phi $ and $\tilde \Psi $:
\begin{equation}
  \label{eq:c34}
\| C_\Psi f\|_{\ell ^p_w} \asymp \| C_{\tilde \Psi } f\|_{\ell ^p_w} \asymp \|
C_\Phi f\|_{\ell ^p_w} \asymp \| C_{\tilde \Phi } f\|_{\ell ^p_w}    
\end{equation}
for $f\in \Hil ^{00}$. This is all that is needed to define
unambiguously a coorbit space $\Hil ^p_w$. 

If \eqref{eq:c34} holds, then all statements of this section,
specifically Theorem~\ref{sec:matbyfram1} and Propositions~\ref{sec:propmatropfram1} --
\ref{prrr} remain true.

One of our main points is that 
these norm equivalences \eqref{eq:c34} always hold for $\mathcal{A}$-localized
frames, as we have seen in 
Section~\ref{sec:localfr0}. In addition, \eqref{eq:c34} also hold  for
wavelet frames (with sufficiently many vanishing moments and sufficient
decay) with  the Besov spaces and
Sobolev spaces as  the corresponding coorbit spaces. In fact, one of
the main motivations for wavelets was the investigation of singular
integral operators, see~\cite{FJ90,meyer1}.  However, for wavelet frames the norm equivalences \eqref{eq:c34} require
different arguments that are not covered by our theory of localized
frames.

\subsection{Invertibility}

For the invertibility we can show, as in the Hilbert space setting \cite{xxlrieck11}:
\begin{lemma}  \label{sec:invertmatrrep1}
Let $\Phi$ and $\Psi$ be $\AAa$-localized frames for $\Hil$ satisfying
$\Phi \sim _\mathcal{A} \Psi $
 and $1\leq p\leq \infty$. 
 Let $O : \Hil_{w_1}^{p} \rightarrow \Hil_{w_2}^{p}$ be a bounded, linear operator.

Then  $O$ is bijective, if and only if $ {\mathcal M}_{(\Phi , \Psi)} \left( O \right)$ is bijective as operator from $\range{C_{\Psi}^{p,w_1}}$ to $\range{C_{\Phi}^{p,w_2}}$. 

In this case the matrix associated to the inverse is given by 
    $${\left( {\mathcal M}_{(\Phi , \Psi)} \left( O \right) \right)}^{\dagger} := {\left({\mathcal M}_{(\Phi , \Psi)} \left( O \right)_{|_{\range{C_{\Psi}}}}\right)}^{-1} = {{\mathcal M}_{(\widetilde \Psi , \widetilde \Phi)}} \left( O^{-1} \right). $$
\end{lemma}
\begin{proof}
By Theorem \ref{sec:locopreint1},  $C_{\Phi}^{p,w}$ is a bijection from
$\HHmp$ onto $\range{C_{\Phi}^{p,w}}$, and  $D_{\Psi}^{p,w}$a bijection
from $\range{C_{\Psi}^{p,w}}$ onto $\HHmp$, where 
$w = w_1$ or $w_2$. 

Therefore $O$ is bijective if and only if ${\mathcal M}_{(\Phi ,
  \Psi)} \left( O \right)$ is bijective from
$\range{C_{\Psi}^{p_1,w_1}}$ to $\range{C_{\Psi}^{p_2,w_2}}$. 

Furthermore
$${{\mathcal M}_{(\widetilde \Psi , \widetilde \Phi)}
\left( O^{-1} \right)} {\mathcal M}_{(\Phi , \Psi)}
 \left( O \right) = 
C^{p_1,w_1}_{\widetilde \Psi} \circ O^{-1}  \circ D^{p_2,w_2}_{\widetilde \Phi} C^{p_2,w_2}_{\Phi} \circ O \circ D^{p_1,w_1}_{\Psi} = $$
$$ = C^{p_1,w_1}_{\widetilde \Psi} \circ D^{p_1,w_1}_{\Psi} = G_{\widetilde \Psi, \Psi}, $$
and therefore the projection on $\range{C_{\Psi}^{p_1,w_1}}$.
\ep
\begin{remark}
Note that for $p_1 \neq p_2$, there does not exist a bijective operator
$O : \Hil_{w_1}^{p_1} \rightarrow \Hil_{w_2}^{p_2}$ by Theorem \ref{sec:equivclass1}. 
 \end{remark}

The condition number of a matrix (or an operator) plays an important role in numerical analysis  \cite{luen} and is defined by $\kappa\left(M\right) = \norm{Op}{M} \cdot \norm{op}{M^{-1}}$. For matrices with non-zero kernel we can define the generalized condition number \cite{chewei04} by $\kappa^\dagger\left(M\right) = \norm{Op}{M} \cdot \norm{op}{M^\dagger}$.  By using Lemma  \ref{sec:invertmatrrep1} and Theorem \ref{sec:matbyfram1} it is straightforward to show
$$ \kappa^\dagger \left({{\mathcal M}_{(\Phi ,   \Psi)} \left( O \right)}\right) = \kappa^\dagger\left({G_{\Phi, \Psi}}\right) \cdot \kappa^\dagger\left(G_{\widetilde \Psi, \Psi}\right) \cdot \kappa^\dagger\left( O \right). $$

 \begin{theorem}
  Let $\Psi $ be an $\mathcal{A}$-localized frame for $\Hil$ and $w$
  an $\mathcal{A}$-admissible weight. Assume that $O: \Hil \to \Hil $
  is invertible and that $ \mathcal{M}_{(\Psi  , \tilde \Psi)} ( O )
  \in \mathcal{A}$. Then $O$ is invertible simultaneously on all
  coorbit spaces $\Hil ^p_w$, $1\leq p \leq \infty $. 
 \end{theorem}

 \begin{proof}
By Lemma~\ref{sec:invertmatrrep1} the matrix of $O^{-1}$ is given by    
    $\mathcal{M}_{(\widetilde \Psi ,  \Psi)} ( O^{-1}) = \Big(
  {\mathcal M}_{(\Psi , \tilde \Psi)} ( O ) \Big)^{\dagger} $. Since
$\mathcal{A}$ is closed with respect to taking a pseudo-inverse and
${\mathcal M}_{(\Psi , \tilde \Psi)} ( O )\in \mathcal{A}$, it follows
that also  $\mathcal{M}_{(\widetilde \Psi ,  \Psi)} ( O^{-1}) \in
\mathcal{A}\subseteq \BL{\ell_{w}^{p},\ell_{w}^{p}}$. By
Lemma~\ref{sec:boundedopequiv1} $O^{-1}$ is therefore bounded on $\Hil
^p_w$ for $1\leq p\leq \infty $. 
 \ep

This result is special for intrinsically localized frames and fails
for wavelet frames. 

\begin{figure}[ht]
\center
	\begin{picture}(100,60)
	\put(10,50){$\Hil_{w_1}^{p_1}$}
	\put(90,50){$\Hil_{w_2}^{p_2}$}
	\put(18,51){\vector(4,0){70}}
	\put(50,52){$O$} 
	\put(13,47){\vector(0,-4){30}}
	\put(11,17){\vector(0,4){30}}
	\put(93,47){\vector(0,-4){30}}
	\put(91,17){\vector(0,4){30}}
	\put(0,10){$\range{C_{\Psi}}$} 
	\put(0,30){$D_{\Psi}$}
	\put(15,30){$C_{\Psi}$} 
	\put(80,30){$D_{\Phi}$} 
	\put(95,30){$C_{\Phi}$} 
	\put(82,10){$\range{C_{\Psi}}$}
	\put(20,11){\vector(4,0){60}}
	\put(38,5){$M = \mathcal M_{(\Psi,\Phi)}(O)$}
	\end{picture}
\caption{All operators in the diagram are bijective, if $M$ or
  equivalently $O$, is bijective.} \label{fig:matop1} 
\end{figure}
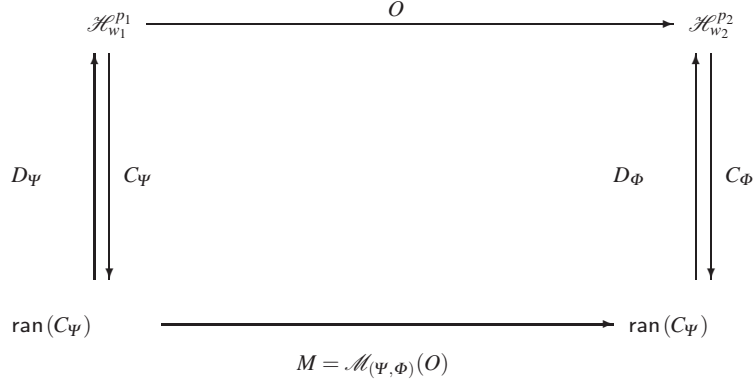

\section{Outlook}  \label{sec:soopeq0}

This manuscript was motivated by many discussions of the first author with applied
scientists who  work  on the numerical solution of  integral
equations. In applications in acoustics,   the solutions of the
Helmholtz equation are of particular  importance, see
e.g. \cite{KasessAPAC2015},  they are  used for example  for the numerical
estimation of head-related transfer functions \cite{Kreuzeretal09,ziegelwangeretalJASA2015}. 

In general the problem of solving an integral equation  can be seen as
solving a linear equation
\begin{equation} \label{sec:solvopeq1}
~ O \cdot f = g ~, 
\end{equation}
where the  operator $O$ models the physical system and the right-hand
side is given and the solution $f$ is to be determined.   For the important  example of   sound propagation, 
the right-hand side $g$ is often called the  load vector. 
It is usually assumed  that $f,g,$ are in some appropriate  function spaces. 

For the numerical treatment  of such operator equations one needs a
reduction to a discrete formulation. This is often done with a
so-called  Galerkin scheme \cite{sauter2010boundary}.  As a first step, either the
boundary of the considered space or the whole space itself are
separated in finite patches or finite volume elements. This leads to
the Boundary Element Method \cite{Gauletal03} or the  Finite Element Method
\cite{brennscott1}.
In the Galerkin scheme one first finds  the matrix
$M=\mathcal{M}_{(\Psi , \Psi )}(0)$ corresponding to
the operator $O$ with respect to a given basis or frame $\Psi
$. Instead of  solving the operator equation $Of=g$, one converts
\eqref{sec:solvopeq1} into a matrix equation as follows. 
$$ O f = g \Longleftrightarrow \sum \limits_l \left< f, \tilde
  \phi_l\right> O  \phi_l = g  \Longleftrightarrow
	\sum \limits_l
\left< f, \tilde \phi_l\right> \left< O \phi_l , \phi_k \right> =
\left< g , \phi_k \right>  \Longleftrightarrow  $$ 
\begin{equation} \label{sec:matrixeq1} 
{\mathcal M}_{(\Phi , {\Phi})} \left( O \right) \cdot C_{\widetilde \Phi} f = C_\Phi g.
\end{equation}
Here $M = {\mathcal M}^{(\Phi , \tilde{\Phi})}
\left( O \right)$ is  called the \emph{system matrix} or
\emph{stiffness matrix}.   

In the  setting of localized frames, the natural  function spaces  are
the Banach spaces $\HHmp$. By the results above  this is equivalent to
the vectors in the matrix equation \eqref{sec:matrixeq1} to be contained
in some  $\ell_w^p$-space. 

In finite and boundary element approach the system $\Phi$ is usually a spline-like
basis \cite{Gauletal03}. Recently wavelet bases \cite{DahSch99a}, but
also frames have been applied, e.g. in \cite{Stevenson03,harbr08}.  
Currently the potential of other  frames is investigated  for solving operator equations  in
acoustics, such as   $\alpha$-modulation
frames \cite{BaySpe2014}.
Note, however,  that neither wavelet frames  nor $\alpha$-modulation
frames are  localized in the sense of Definition~\ref{defloc}.  As mentioned above, as long as \eqref{eq:c34} is fulfilled most results of Section \ref{sec:matreplocalfr0} can still
be applied. 

To use a numerical solver, it is necessary to perform a further
reduction to a finite-dimensional matrix equation. This means that we  have to find a good
finite dimensional approximation  of $\mathcal{M}_{(\Phi , \Phi )}(O)$.  
This is done by restricting $\mathcal{M}_{(\Phi , \Phi )}(O)$ to a
suitable finite-dimensional subspace. Specifically,  let $\left\{ P_n \right\}$
be a bounded sequence of finite-rank orthogonal projections in $\BL{\Hil}$ with the property
that $P_n x \rightarrow x$ for all $x \in \Hil$ and $n \rightarrow \infty$. 
Assume that  $A \in \BL{\Hil}$
is invertible. Consider  
\begin{equation} \label{sec:finitsecdef1} 
P_n A P_n x = P_n y
\end{equation}
and solve for $x_n = \left( P_n A P_n  \right)^{-1} P_n y$. This is the classical 
 \emph{projection method}~\cite{gohberg1}.  The projection method for
 $A$  is said to \emph{converge}, if 
 for all $y \in \Hil$ there exists a unique solution $x_n$
 to~\eqref{sec:finitsecdef1} with $x_n \rightarrow A^{-1} y$. This is 
 the case \cite{gohberg1} if and only if the matrices $A_n = P_n A P_n$ have
 uniformly bounded inverses, i.e. $\sup \limits_{n \ge N} \norm{}{A_n^{-1}}
 < \infty$. In particular, if 
 $\norm{}{I - A} < 1$, then the method converges. The special case 
 when  $P_n$ is the
 orthogonal projection on the first $n$ coordinates in $\ell^2$ 
 is called the \emph{finite section} method.
 In numerical analysis, this approximation  scheme is often called the
 Galerkin scheme. 

If   $\left| \left< A x, x \right> \right| \ge c \cdot \norm{}{x}^2$,
then $A$ is invertible, $\norm{}{A^{-1}} \le c^{-1}$ and it is easy to
see that the
projection method converges. 

A convergence analysis of the finite section method in weighted $\ell
^p$-spaces is carried out in~\cite{grrzst10}. The methods are closely
related to the methods used for the analysis of localized frames.

The projection method can be combined with frames in several ways:
\begin{itemize}
\item The ``naive approach'': assume that $K= \mathbb{Z}$ and  choose  ${M_N}_{k,l} =
  \mathcal{M}_{(\Phi , \Phi )}(O)_{k,l}$ for $|k|,|l|\le N$. This
  corresponds to a finite section method. 
\item Subspace selection: Choose a sequence  $K_N$ of finite  subsets of $K$ with the following properties: 
\begin{enumerate}
\item[(i)] $K_i \subseteq K_j$ for $i \le j$ and $\bigcup _{i=1}^\infty K_i
  =  K$. 
\item[(ii)] The space $V_N := \mathrm{span}\, \{ \psi_k | k  \in K_N  \}$ has dimension $N$.
\end{enumerate} 
For a frame $\Psi $ it may happen that $\mathrm{dim} \, V_N <
\mathrm{card} (K_N)$, but the set $\Psi^{(N)} := \{ \psi_k | k  \in K_N  \}$ is
always a frame for $V_N$. (This need not be the case when
$\mathrm{dim} \,  V_N = \infty $).   
For the numerical treatment the condition numbers of the transforms, i.e. the quotients of the frame bounds, has to be controlled. 
Therefore we consider $\Psi^{(N)} = \left\{ \psi_k \left| k \in K_N
  \right. \right\}$ being a frame for $V_N$ with bounds $C$,$D$ independent of $N$. 
	This is called a subframe in \cite{harbr08}.
	We denote the canonical dual on $V_N$ by
$\widetilde{\Psi}^{(N)} := \left\{ \tilde{\psi}_k^{(N)} \right\}$. 
Then use the projection $P_N f = \sum \limits_{k \in K_N} \left< f ,
  \psi_k \right> \tilde{\psi}_k^{(N)} = \sum \limits_{k \in K_N}
\left< f , \tilde{\psi}_k^{(N)} \right> \psi_k$ and solve
\eqref{sec:finitsecdef1}. Since  $\range{{\mathcal
    M}^{(\Phi , {\Phi})} \left( O \right) } \subseteq
\range{C_{\Phi}}$,  this set-up   leads  exactly to the formulation
in  \eqref{sec:matrixeq1}. 

In concrete applications it is a non-trivial problem to  find  index
sets such that the approximation method converges and at the same time  is numerically efficient.
For wavelet frames this problem can be tackled with  a
multi-resolution  approach  with a basis property on each scale,  see,
e.g., \cite{harbr08}. 
\end{itemize}
We note that the matrix $M_N$ cannot have full rank  whenever   the
frame $\{\psi _k | k\in K_N\}$  is redundant for $V_N$.   
By Lemma \ref{sec:invertmatrrep1} the equation \eqref{sec:matrixeq1}
still has a unique solution, although the matrix is not
invertible. For the efficient solution of \eqref{sec:matrixeq1}, even  for
frames,  one can apply Krylov subspace
methods, such as the conjugate gradient method \cite{harbr08}.  
Other possible methods include versions of Richardson iterations
\cite{dafora05} or steepest descent methods \cite{Dahlkeetal07c}. \\

\begin{acknowledgement}
The first author was in part supported by the 
START project FLAME Y551-N13
of the Austrian Science Fund (FWF) and  the DACH project BIOTOP
I-1018-N25 of Austrian Science Fund (FWF). The second author
acknowledges the support of the  FWF-project P 26273-N25. P.B.\  wishes to thank
NuHAG for the hospitality as well as the availability of its webpage.  
He also thanks Dominik Bayer, Gilles Chardon, Stephan Dahlke, Helmut
Harbrecht, Wolfgang Kreuzer, Michael Speckbacher and Diana Stoeva for related interesting discussions. 
\end{acknowledgement}

{
\small


}
\end{document}